\documentclass[11pt,reqno]{article}

\usepackage{enumitem,kantlipsum}
\usepackage[table]{xcolor} 
\usepackage[table]{xcolor} 
\usepackage{booktabs} 
\usepackage{adjustbox} 
\usepackage{cancel}
\usepackage{epigraph}
\usepackage{geometry}                
\usepackage{subcaption}
\setlist[enumerate, 1]{label=\arabic*.  }

\usepackage[dvipsnames]{xcolor}
\usepackage{graphicx}
\usepackage{amssymb}
\usepackage{amsthm}
\usepackage{epstopdf}
\RequirePackage{doi}
\usepackage{hyperref}
\usepackage{graphicx}
\usepackage{forest}
\usepackage{tikz}
\usepackage{tikz-cd}
\usepackage{stmaryrd}
\usepackage{enumitem}
\usepackage{comment} 
\usepackage{newtxmath}
\usepackage{nicefrac}
\usepackage{pgfplots}
\pgfplotsset{compat=1.18}

 \usepackage{pdfpages}

\definecolor{USred}{cmyk}{0,1.00,0.65,0.34}
\definecolor{darkblue}{rgb}{0.0627, 0.1539, 0.4510}
\definecolor{darkgreen}{rgb}{0.0627,0.4510,0.1539}

\definecolor{darkorange}{cmyk}{0,0.55,1,0.20}

\renewcommand{\emph}[1]{{\textcolor{USred}{\em #1}}}

\hypersetup{
  colorlinks,
  citecolor=USred,
  linkcolor=USred,
  urlcolor=USred}

\def\SS{\mathbb{S}}

\DeclareMathOperator{\blob}{RecordCode}

\DeclareMathOperator{\rec}{rec}

\DeclareMathOperator{\mrec}{mrec}

\newcommand{\planted}{R_{\bullet}}

\theoremstyle{definition}
\newtheorem{de}{Definition}[section]
\theoremstyle{plain}
\newtheorem{thm}[de]{Theorem}

\newtheorem{lem}[de]{Lemma}

\newtheorem{p}[de]{Proposition}
\newtheorem{cor}[de]{Corollary}
\newtheorem{ex}[de]{Example}

\theoremstyle{definition}

\title{The genesis sequence, tree records and endofunctions.}
\author{Enrica Duchi, Adrián Lillo, Pablo Puerto, Mercedes Rosas, Stefan Trandafir}

\begin{document}

\maketitle

\begin{abstract}
In this work, we present a series of bijections that reveal the deep connections between the concepts of  tree records, the girth of a connected endofunction, and the genesis sequence, the first sequence in the OEIS. We use these results to derive the generating functions for the tree and forest record numbers, expressing them in terms of the Cayley's tree function.  Finally, we provide a new proof for Cayley's forest formula.
\end{abstract}

\vspace{-.2cm}
\epigraph{\em
        This paper is dedicated to Neil J. A. Sloane, renowned number-crusher, for his creation and tireless curation of the Online Encyclopedia of Integer Sequences.
}

\tikzset{
    record/.style={
        circle,
        draw=black,
        fill=black,
        inner sep=1.5pt,
        label={#1},
        solid
    },
    record/.default={},
    non-record/.style={
        circle,
        draw = black,
        fill = white,
        inner sep=1.5pt,
        solid
    },
    every label/.style={
            font=\tiny
        }
}


\section*{Introduction}

Trees and forests are among the most fundamental structures in discrete mathematics. 
This article builds on a beautiful connection between tree records, heights of nodes of trees, and endofunctions discovered using Sloan's online encyclopedia of integer sequences, and that goes back to the very first sequence of the encyclopedia, the ``genesis sequence".

In the Numberphile podcast episode “The Number Collector” \cite{podcast_Sloan}, Neil Sloane recounts how he encountered the genesis sequence of the OEIS, the sequence  ``Normalized total height of all nodes in all rooted trees with $n$ labeled nodes",  during his PhD research, and  how it inspired the creation of the OEIS.
A short version of this story appears at the encyclopedia's wiki that we quote textually:
``The sequence database was begun by Neil J. A. Sloane  in early 1964 when he was a graduate student at Cornell University in Ithaca, NY. He had encountered a sequence of numbers while working on his dissertation, namely 1, 8, 78, 944, ... (now entry \url{https://oeis.org/A000435} in the OEIS), and was looking for a formula for the $n$-th term, in order to determine the rate of growth of the terms.
He noticed that although several books in the Cornell library contained sequences somewhat similar to this, this particular sequence was not mentioned. In order to keep track of the sequences in these books, NJAS started recording them on file cards, which he sorted into lexicographic order."
(From \url{https://oeis.org/wiki/Welcome#OEIS:_Brief_History})

In this work, we build on our study of records  of rooted trees and rooted forests initiated in \cite{LRT-GF-records, LRT-Weary}.  Some quick, necessary definitions to be detailed later. A \emph{tree record} of a rooted labeled tree is a node whose label is greater than the label of all other nodes  on the path to the root. The \emph{height} of a  node is the number of edges it has to cross to get there.

An  \emph{endofunction} is a map from some $[n] $ to itself. Its functional graph decomposes into connected components, each consisting of a directed cycle of rooted trees.  The \emph{girth} of a connected endofunction is the number of edges in the unique cycle of its functional graph.

The material presented in this work began with the discovery that the genesis sequence also counts the total number of records in all rooted trees with 
$n$ nodes, a fact that we establish in Theorem \ref{thm:catalysts}. On the other hand, we wanted to provide combinatorial proofs for some interesting expressions relating tree records to the Cayley tree function presented in \cite{LRT-GF-records}. This is the content of the second part of the work, and follows from the combinatorial framework developed to deal with the notions of tree records, nodes at a given height, and the girth of endofunctions.

In Theorem \ref{thm:connected endofunctions} we present a bijection, inspired on Joyal's proof of Cayley's formula, between  the set of rooted trees of order $n$ with a distinguished record at height $k-1$ and the set of connected endofunctions on $[n]$ of girth  $k$.  Figure  \ref{fig:endomorphism bijection} illustrates this bijection. The existence of this bijection implies that the total number of records in all rooted trees with $n$ nodes is equal to the
number of connected endofunctions on $[n].$

In Theorem   \ref{Thm:endo_girth}  we present a bijection  between  connected endofunctions on $[n]$ of girth $m$ with at least $k + 1$ records and connected endofunctions on $[n]$ of girth $m + 1$ with at least $k$ records. 
In particular, this result implies that the number of endofunctions on $[n]$ with $k$ records and girth $m$ only depends on $m+k$. Moreover, combining Theorem \ref{thm:connected endofunctions} and Theorem   \ref{Thm:endo_girth} we conclude that the set of rooted trees of order $n$ with a distinguished record at height $k-1$
is in bijection with the set of trees of order $n$ with at least $k$ records.
From Theorem   \ref{Thm:endo_girth} we also derive a direct combinatorial proof for  the record number formula of \cite{LRT-GF-records} for the number of  rooted   trees labeled with  $[n]$ and  with $k$ records:
    \[
    \planted(n,k) = 
       k \,(n-1)\cdots
    (n-k+1) \, n^{n-k-1}.
     \] 

In Section  \ref{se:genesis} we focus on the connection between tree records, connected endofunctions, and the genesis sequence of the encyclopedia. We present a bijection that shows that the total number of non-root records in all rooted trees labeled with $[n]$ is
equal to the normalized total height of all nodes in all rooted trees with $n$ labeled nodes, an expression that defines the genesis sequence of the OEIS. This is done in  Theorem \ref{thm:catalysts}.

In Section \ref{se:record_code_endofunction} we study the endofunction defined by the record code of a rooted tree.
First, we recall the notions of bonsai decomposition and record code of a rooted tree. Then, we consider the relation between record distances of a rooted tree $T$ and the order structure defined by $T$. In particular, we find a closed formula for the number of rooted trees where we fix the distances between records.

In Section \ref{se:record_gf}, we translate the combinatorial framework developed during this exploration to the language of  generating functions. Let  $\mathcal{T}(z)$ be the ubiquitous Cayley tree function, and let  $\mathcal{R}_{\ge k}(z)$ be the generating function for rooted trees with at least $k$ records weighted by the number of nodes, and $\mathcal{R}_{k}'(z)$ for trees with $k$ records and a distinguished node. Then, we show that
\begin{align*}
\mathcal{R}_{\ge k}(z) = \frac{\mathcal{T}^k(z)}{k},
\qquad \text{ and } \qquad 
\mathcal{R}_{k}'(z) = \frac{\mathcal{T}^k(z)}{z}.
\end{align*}
and derive several related identities from these two expressions, which were first obtained in \cite{LRT-GF-records} through manipulations of generating functions.

Section \ref{App:Cayley_forest} presents a coda for our work. In 1889, Cayley~\cite{Cayley} discovered that the number of forests on $n$ labeled vertices with roots specified by $[k]$ is $k\, n^{\,n-k-1}$. This result is known as \emph{Cayley’s forest formula} and has many definite proofs.  {\it Proofs from THE BOOK}~\cite{thebook} presents two clean arguments for this identity: an inductive proof by Riordan and Rényi, and a double-counting proof by Pitman. Peter Shor gives another proof for this formula  in \cite{Shor}, emphasizing tree inversions. We dare to propose yet another proof, this time based on the notion of tree records. Cayley's forest formula is relevant to our work because  our proof of the tree record number rests on it.  In particular, and for fixed $k$, Lemma \ref{lem:blob-iff} characterizes the record codes for trees where the children of the root are labelled with $[k]$ in terms of the lexicographic order.


\section{Connected endofunctions and tree records.}
We begin with some basic conventions and definitions followed in this work.
 We write $[n]$ for the set $\{1, 2, \ldots, n\}$ and $[n]_0$ for $\{0, 1, \ldots, n\}$.
    As all structures we consider are labeled, we simply use the terms tree and forest to refer to labeled trees and labeled forests.  

 Let $T$ be a rooted tree. The number of nodes in a tree or a forest (not counting $\circ$ whenever  it happens to appear as one of the labels) defines its \emph{order}.
A node $v$ of $T$ is a \emph{record} if it has the largest label along the unique path from $v$ to the root. Otherwise we say that $v$  is a non-record.  
The \emph{tree record number} $\planted(n,k)$ is defined as the number of rooted trees of order $n$ with precisely $k$ records.
The \emph{height} of a node of $T$ is the length of the path (number of edges) from it to the root. 
The \emph{total height} of $T$ is the sum of the heights of the nodes of $T$.

An \emph{endofunction} on $[n]$ is a function $f : [n] \to [n]$ for some natural number $n$. Its \emph{functional graph} is  the directed graph with nodes $[n]$, and one directed edge $(i,f(i))$ for each $i$ in $[n]$.
Given a rooted tree $T$, we define the \emph{parent map} of $T$ as the map that sends each non-root vertex to its parent, and the root to itself.  Observe that this is a natural way of embedding the set of rooted trees in the set of endofunctions.
An endofunction  is said to be \emph{connected} if its functional graph is connected after ignoring orientations. The \emph{girth} of a connected endofunction is the number of edges in the unique cycle of its functional graph.

\subsection{The height of a record and the girth of an endofunction}

Through the lens of Joyal’s bijection \cite{Joyal}, we reveal an elegant relationship between record numbers and connected endofunctions that constitutes the starting point of our journey. Fix two nonnegative integers $k \le n$.

\begin{thm} 
\label{thm:connected endofunctions}
The number of rooted trees of order $n$ with a distinguished record at height $k-1$ coincides with the number of connected endofunctions on $[n]$ of girth $k$.
\end{thm}

\begin{figure}[h!]
    \centering

\resizebox{0.75\textwidth}{!}{$
  \vcenter{\hbox{\usetikzlibrary{positioning}
\tikzset{
  circ/.style={
    circle, draw, semithick, minimum size=1.5mm, inner sep=0pt, fill=black,
    }
}

\begin{tikzpicture}[
    arr/.style={semithick}
]
    \newcommand{\len}{1}

    \node[circ, label={right:$\textcolor{USred}{5}$}, fill=USred,USred] (s1) at (0,0) {};
    \node[circ, label={right:$\textcolor{USred}{7}$}, fill=USred,USred] (s2) at ($(s1) + (-90:\len)$) {};
    \node[circ, label={right:$\textcolor{USred}{3}$}, fill=USred,USred] (s3) at ($(s2) + (-90:\len)$) {};
    \node[circ, label={right:$\textcolor{USred}{9}$}, fill=USred,USred] (s4) at ($(s3) + (-90:\len)$) {};

    \node[circ, label={right:$16$}] (s5) at ($(s4) + (-90:\len)$) {};
    
    \node[circ, label={right:$11$}] (s6) at ($(s3) + (-45:\len)$) {};
    \node[circ, label={right:$6$}] (s7) at ($(s6) + (-45:\len)$) {};
    \node[circ, label={right:$4$}] (s8) at ($(s6) + (-90:\len)$) {};
    \node[circ, label={right:$15$}, label distance=0] (s16) at ($(s3) + (-135:\len)$) {};

    \node[circ, label={right:$1$}] (s10) at ($(s2) + (-135:\len)$) {};
    \node[circ, label={right:$14$}] (s11) at ($(s10) + (-135:\len)$) {};
    \node[circ, label={above:$8$}] (s12) at ($(s10) + (-180:\len)$) {};
    \node[circ, label={right:$2$}] (s9) at ($(s2) + (-45:\len)$) {};

    \node[circ, label={right:$13$}] (s13) at ($(s1) + (-135:\len)$) {};
    \node[circ, label={right:$10$}] (s14) at ($(s1) + (-45:\len)$) {};
    \node[circ, label={right:$12$}] (s15) at ($(s14) + (-45:\len)$) {};

    \begin{scope}[USred]
    \draw[arr] (s2) -- (s1);
    \draw[arr] (s3) -- (s2);
    \draw[arr] (s4) -- (s3);
    \end{scope}

    \draw[arr] (s5) -- (s4);
    
    \draw[arr] (s6) -- (s3);
    \draw[arr] (s7) -- (s6);
    \draw[arr] (s8) -- (s6);
    \draw[arr] (s16) -- (s3);
    
    \draw[arr] (s10) -- (s2);
    \draw[arr] (s11) -- (s10);
    \draw[arr] (s12) -- (s10);
    
    \draw[arr] (s13) -- (s1);
    \draw[arr] (s14) -- (s1);
    \draw[arr] (s15) -- (s14);
    \draw[arr] (s9) -- (s14);

    \draw[draw=USred, semithick, fill=white] (s4) circle [radius=1.6mm];
    \node[circ, fill=USred, draw=USred] (s4bis) at (s4) {};

\end{tikzpicture}}}%
  \hspace{1cm}%
  \longleftrightarrow%
  \hspace{1cm}%
  \vcenter{\hbox{\begin{tikzpicture}[
    circ/.style={circle, draw, semithick, minimum size=1.5mm, inner sep=0pt, fill=black},
    arr/.style={-stealth, semithick}
]
    \newcommand{\len}{0.95}
    \begin{scope}[USred]    
    \node[circ,  label={right:$\scriptsize 7$}, fill=USred,USred] (c1) at (0, 0) {};
    \node[circ, label={right:$\scriptsize 5$}, fill=USred,USred] (c2) at ($(c1) + (-135:\len)$) {};
    \node[circ,  label={right:$\scriptsize 9$}, fill=USred,USred] (c3) at ($(c2) + (-45:\len)$) {};
    \node[circ, label={right:$\scriptsize 3$}, fill=USred,USred] (c4) at ($(c3) + (45:\len)$) {};
    \end{scope}

    \node[circ, label={right:$\scriptsize 1$}] (c5) at ($(c1) + (90:\len)$) {};
    \node[circ, label={right:$\scriptsize 8$},] (c6) at ($(c5) + (45:\len)$) {};
    \node[circ, label={right:$\scriptsize 14$},] (c7) at ($(c5) + (135:\len)$) {};

    \node[circ, label={right:$\scriptsize 16$},] (c8) at ($(c3) + (-90:\len)$) {};

    \node[circ, label={right:$\scriptsize 15$},] (c9) at ($(c4) + (45:\len)$) {};
    \node[circ, label={above right:$\scriptsize 11$},] (c10) at ($(c4) + (-45:\len)$) {};
    \node[circ, label={right:$\scriptsize 6$},] (c11) at ($(c10) + (0:\len)$) {};
    \node[circ, label={right:$\scriptsize 4$},] (c12) at ($(c10) + (-90:\len)$) {};
    
    \node[circ, label={right:$\scriptsize 10$},] (c13) at ($(c2) + (135:\len)$) {};
    \node[circ, label={right:$\scriptsize 13$},] (c14) at ($(c2) + (-135:\len)$) {};
    \node[circ, label={right:$\scriptsize 2$},] (c15) at ($(c13) + (135:\len)$) {};
    \node[circ, label={right:$\scriptsize 12$},] (c16) at ($(c13) + (-135:\len)$) {};

    \begin{scope}[USred]
    \draw[arr] (c1) -- (c2);
    \draw[arr] (c2) -- (c3);
    \draw[arr] (c3) -- (c4);
    \draw[arr] (c4) -- (c1);
    \end{scope}

    \draw[arr] (c5) -- (c1);
    \draw[arr] (c6) -- (c5);
    \draw[arr] (c7) -- (c5);
    
    \draw[arr] (c16) -- (c13);
    \draw[arr] (c14) -- (c2);
    \draw[arr] (c15) -- (c13);
    \draw[arr] (c13) -- (c2);

    \draw[arr] (c8) -- (c3);

    \draw[arr] (c9) -- (c4);
    \draw[arr] (c10) -- (c4);
    \draw[arr] (c11) -- (c10);
    \draw[arr] (c12) -- (c10);
    
\end{tikzpicture}}}
  $}
\caption{The bijection of Theorem \ref{thm:connected endofunctions}. The distinguished record $v = 9$ is circled. The path from $v$ to the root (on the left) and the only cycle in the resulting endofunction (on the right) are highlighted.}
    \label{fig:endomorphism bijection}
\end{figure}
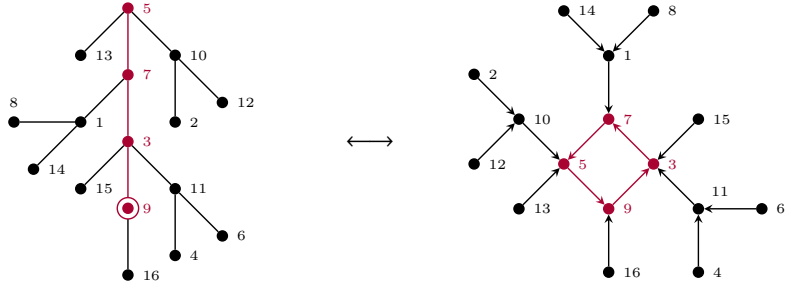

\begin{proof} We present an explicit bijection to prove this fact. Let $T$ be a rooted   tree with a distinguished record $v$,  root $r$. Let $p$ be its parent map.
We send the pair $(T, v)$ to the endofunction $f : [n] \to [n]$ defined by 
\[
f(x) = p(x) \quad \text{if} \quad x\neq r, \qquad f(r) = v.
\]
This is, $f$  sends the root $r$ to the distinguished record $v$, and any other vertex to its parent in $T$.  See Figure \ref{fig:endomorphism bijection}.
Note that $f$ is a connected endofunction, as
$T$ is a spanning tree of the graph of $f$.

We show that the map  $(v, T) \mapsto f$ is a bijection by exhibiting the inverse map. 
Let $f$ be a connected endofunction of $[n]$, and let $\sigma =(x_1 \ x_2 \ \cdots \ x_m)$ be the only cycle of $f$. First, we need to determine the vertices $v$ and $r$. These are given by
$
v = \max_i x_i,$ and $r = \sigma^{-1}(v).
$
 We then reconstruct the tree $T$ by first removing the edge that connects $v$ with $r$ in the graph of $f$, and then rooting the resulting tree at $r$. Note that $v$ is a record of $T$: the vertices in the path from $v$ towards the root are precisely the elements of $\sigma$, and we have defined $v$ as the maximum of these vertices. Consequently, the height of $v$ is exactly one less than the girth of $f$.

\qedhere
\end{proof}

\begin{cor} 
\label{cor:connected endofunctions}
The total number of records in all rooted trees with $n$ nodes is equal to the number of connected endofunctions on $[n]$. 
\end{cor}

The sequence counting connected endofunctions has been studied in depth, see \cite{Flajolet} and \url{https://oeis.org/A001865} for more information. The refined sequence that keeps track of the endofunction's girth also appears in the encyclopedia \url{https://oeis.org/A201685}. Theorem \ref{thm:connected endofunctions} and Corollary \ref{cor:connected endofunctions} present two new combinatorial interpretations for these sequences of numbers based on the notion of tree record.

\subsection{Some consequences of our bijections}

Let $f$ be an endofunction on $[n]$.   The set of elements of the form $f^k(i)$ for some $k\ge0$ is referred to as the \emph{orbit} of  $i$, for any fixed element $i$ of $[n]$. If the endofunction $f$ is connected we refer to the \emph{cycle} of $f$ as the cyclic permutation obtained by restricting $f$ to the elements of its unique cycle.
Finally, an element $i \in [n]$ is said to be a \emph{record} if $ i$ is bigger or equal that all the elements of its orbit. Fix three nonnegative integers $k,m \le n$.

\begin{thm}
\label{Thm:endo_girth}
    The number of connected endofunctions on $[n]$ of girth $m$ and at least $k + 1$ records coincides with the number of connected endofunctions on $[n]$ of girth $m + 1$ and at least $k$ records.  
\end{thm}

\begin{proof} Again, we present an explicit bijection to prove this fact.
    Let $f$ be an connected endofunction of girth $m$ and at least $k+1$ records. Let $v$ be the $k$-th greatest record of $f$.  
    Observe that there is exactly one record on the cycle $\sigma$ of $f$: the smallest record of the endofunction, and that it is smaller than the $k$-th greatest record. Therefore, $v$ does not belong to $\sigma$.
    
Let $w$ be the first element in the orbit of $v$ under $f$ such that $f(w)$ belongs to $\sigma$.
Denote $x = \sigma^{-1}(f(w))$.
We define an endofunction $g$ by 
\[
g(z) = f(z) \quad \text{if } z \neq x, \qquad g(x) = w.
\]
It is not hard to verify that $g$ is a connected endofunction of girth $m+1$ and its $k$ greatest records coincide with those of $f$. Reversing the bijection is straightforward by first identifying $v$ and $w$.
\end{proof}

\begin{figure}
\[
\centering
\vcenter{\hbox{
\begin{tikzpicture}[
    circ/.style={circle, draw, semithick, minimum size=1.5mm, inner sep=0pt, fill=black},
    rec/.style={circle, draw=USred,  minimum size=1.5mm, inner sep=0pt, fill=USred},
    arr/.style={-stealth, semithick}
]
    \def\sepFirstLevel{0.8}
    \def\len{0.8}
    \node[rec, label={left:3}] (3) at (0, 0) {};
    \node[rec, label={right:$4$}] (4) at ($(3) + (-50:\len)$) {};
    \node[rec, label={right:$5$}, label={left:\textcolor{USred}{$v$}}, label={above:\textcolor{USred}{$w$}}] (5) at ($(3) + (-130:\len)$) {};
    \node[rec, label={below:$6$}] (6) at ($(5) + (-120:0.9*\len)$) {};
    \node[circ, label={below:$2$}] (2) at ($(5) + (-60:0.9*\len)$) {};
    
    \node[circ, label={below:$1$}] (1) at ($(4) + (-120:0.9*\len)$) {};
    \node[rec, label={below:$7$}] (7) at ($(4) + (-60:0.9*\len)$) {};

\begin{scope}[-stealth]    
   \draw (3) edge[-stealth, out=45, in=135, looseness=12] (3); 
    \draw (6) -- (5);
    \draw (2) -- (5);
    \draw (5) -- (3);
    \draw (4) -- (3);
    \draw (1) -- (4);
    \draw (7) -- (4);
\end{scope}
\end{tikzpicture}}}
\qquad
\vcenter{\hbox{
\begin{tikzpicture}[
    circ/.style={circle, draw, semithick, minimum size=1.5mm, inner sep=0pt, fill=black},
    rec/.style={circle, draw=USred,  minimum size=1.5mm, inner sep=0pt, fill=USred},
    arr/.style={-stealth, semithick}
]
    \def\sepFirstLevel{0.8}
\def\len{0.8}
    \node[circ, label={above:3}] (3) at (0, 0) {};
    \node[circ, label={above:$4$}] (4) at ($(3) + (0:\len)$) {};
    \node[rec, label={above:$5$}] (5) at ($(3) + (-180:\len)$) {};
    \node[rec, label={above:$7$}] (7) at ($(4) + (30:\len)$) {};
    \node[circ, label={above:$1$}] (1) at ($(4) + (-30:\len)$) {};
    \node[rec, label={above:$6$}, label={left:\textcolor{USred}{$v$}},
    label={right:\textcolor{USred}{$w$}}] (6) at ($(5) + (150:\len)$) {};
    \node[circ, label={above:$2$}] (2) at ($(5) + (210:\len)$) {};

    \begin{scope}[-stealth]
    \draw (6) -- (5);
    \draw (2) -- (5);
    \draw (4) -- (3);
    \draw (5) edge[bend right] (3);
    \draw (3) edge[bend right] (5);
    \draw (1) -- (4);
    \draw (7) -- (4); 
    \end{scope}
\end{tikzpicture}
}}
\qquad
\vcenter{\hbox{
\begin{tikzpicture}[
    circ/.style={circle, draw,  minimum size=1.5mm, inner sep=0pt, fill=black},
    rec/.style={circle, draw=USred,  minimum size=1.5mm, inner sep=0pt, fill=USred},
    arr/.style={-stealth}
]
    \def\len{0.8}
    \node[circ, label={below:$\scriptsize 3$},] (6) at (0, 0) {};
    \node[circ, label={below:$\scriptsize 5$},] (5) at ($(6) + (180:\len)$) {};
    \node[rec, label={left:$\scriptsize 6$},] (3) at ($(6) + (120:\len)$) {};
    \node[circ, label={above:$4$}, label={[yshift=2pt, xshift=2pt]below left:\textcolor{USred}{$w$}}] (4) at ($(6) + (-30:\len)$) {};
    \node[rec, label={above:$7$}, label={right:\textcolor{USred}{$v$}}] (7) at ($(4) + (0:\len)$) {};
    \node[circ, label={above:$1$}] (1) at ($(4) + (-60:\len)$) {};
    \node[circ, label={below:$2$}] (2) at ($(5) + (-150:\len)$) {};

    \begin{scope}[-stealth]
        \draw (5) -- (6);
        \draw (3) -- (5);
        \draw (6) -- (3);
        \draw (4) -- (6);
        \draw (1) -- (4);
        \draw (7) -- (4);
        \draw (2) -- (5);
    \end{scope}
\end{tikzpicture}}}
\qquad
\vcenter{\hbox{
\begin{tikzpicture}[
    circ/.style={circle, draw, semithick, minimum size=1.5mm, inner sep=0pt, fill=black},
    rec/.style={circle, draw=USred,  minimum size=1.5mm, inner sep=0pt, fill=USred},
    arr/.style={-stealth}
]
    \newcommand{\len}{0.8}
    \begin{scope}    
    \node[circ,  label={right:$\scriptsize 3$}] (c1) at (0, 0) {};
    \node[rec, label={right:$\scriptsize 6$}] (c2) at ($(c1) + (-135:\len)$) {};
    \node[circ,  label={right:$\scriptsize 5$}] (c3) at ($(c2) + (-45:\len)$) {};
    \node[circ, label={right:$\scriptsize 4$}] (c4) at ($(c3) + (45:\len)$) {};
    \end{scope}

    \node[rec, label={right:$\scriptsize 7$},] (15) at ($(c4) + (45:\len)$) {};
    \node[circ, label={right:$\scriptsize 1$},] (11) at ($(c4) + (-45:\len)$) {};
    
    \node[circ, label={right:$\scriptsize 2$},] (10) at ($(c3) + (-90:\len)$) {};

    \begin{scope}
    \draw[arr] (c1) -- (c2);
    \draw[arr] (c2) -- (c3);
    \draw[arr] (c3) -- (c4);
    \draw[arr] (c4) -- (c1);
    \end{scope}

    \draw[arr] (10) -- (c3);
    \draw[arr] (11) -- (c4);
    \draw[arr] (15) -- (c4);

\end{tikzpicture}
}}
\]
    \caption{The repeated application of the bijection in Theorem \ref{Thm:endo_girth}. The records are highlighted. Note that the leftmost endofunction has at least $4$ records and girth $1$, and the rightmost endofunction has at least $1$ record and girth $4$.}
    \label{fig:placeholder}
\end{figure}
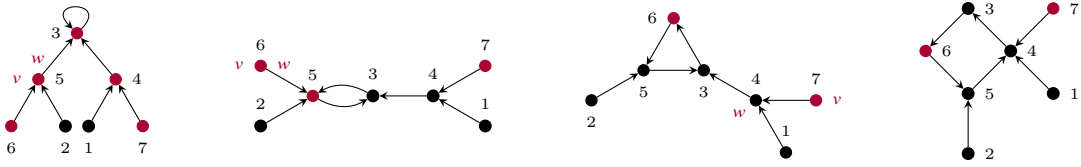

\begin{cor}
\label{pablos_lemma}
The number of connected endofunctions on $[n]$ of girth $m$ and at least $k$ records only depends on $m+k$. In particular, the number of rooted trees with $n$ nodes and at least $k$ records coincides with the number of connected endofunctions on $[n]$ of girth $k$.
\end{cor}
\begin{proof}
Fix two positive integers $k$ and $m$ and let $E_{m,k}$ denote the set of endofunctions of $[n]$ with girth $m$ and at least $k$ records. Applying repeatedly the bijection of Theorem \ref{Thm:endo_girth}, we obtain the following 
bijections:
\[
E_{m, k} \longleftrightarrow E_{m +1, k -1} \longleftrightarrow E_{m+2, k-2} \longleftrightarrow \dots \longleftrightarrow E_{m + k - 1, 1}.
\]
Therefore, the cardinality of $E_{m, k}$ only depends on $m + k$. Thus, with $m=1$, $
|E_{1, k}| = |E_{k, 1}|.$
The result follows by noticing that endofunctions with girth 1 are parent maps of trees, and that any endofunction has at least one record.
\end{proof}

\begin{cor}
The number of rooted trees of order $n$ with at least $k$ records coincides with the number of tree cycles of order $n$ of $k$ rooted trees.
\end{cor}
\begin{proof}
The functional graph of an endofunction on $[n]$ of girth $k$ is a directed cycle of $k$ rooted trees. Therefore, the bijection of Corollary \ref{pablos_lemma} for $m=1$ gives a bijection between rooted trees with at least $k$ records and cycles of $k$ trees.
\end{proof}

\begin{cor}
The number of trees of order $n$ such that the subgraph induced by its records is connected coincides with the number of cycles of order $n$ of $k$ unrooted trees.
\end{cor}
\begin{proof}
Let $T$ be such a tree. Then, for each record $v$ of $T$, the parent node of $v$ is also a record. Applying the previous bijection to $T$, each of the trees of the resulting endofunction has only one record, the root (the element of the cycle). So we can forget these roots to obtain a cycle of unrooted trees. To recover the endofunction it is enough to identify the root each tree, which is its greatest vertex.

In fact, the resulting trees are precisely those of the bonsai decomposition of $T$, see \cite{LRT-Weary} and the coda (Section \ref{App:Cayley_forest}).
\end{proof}

\begin{cor}
The set of rooted trees of order $n$ with a distinguished record at height $k-1$ is in bijection with the set of trees of order $n$ with at least $k$ records. 
\end{cor}
\begin{proof}
    Apply the bijection of Theorem \ref{thm:connected endofunctions} to obtain a connected endofunction on $[n]$ of girth $k$. Then apply the inverse of the bijection of Corollary \ref{pablos_lemma}.
\end{proof}

\begin{cor}
    The set of rooted trees of order $n$ and with a distinguished record $v$ at height $m$ and with at least $k$ records greater than $v$ is in bijection with the set of endofunctions on $[n]$ of girth $m+k+1$.
\end{cor}
\begin{proof}
    Apply the bijection of Theorem \ref{thm:connected endofunctions} to obtain a connected endofunction on $[n]$ of girth $m+1$, noting that it preserves the $k$ records greater than $v$. Then apply the bijection of Corollary \ref{pablos_lemma}.
\end{proof}

\subsection{The tree record formula}

There exist an elegant closed formula for the tree record numbers, first proved with the use of generating functions in \cite{LRT-GF-records}. In this section, we present a direct combinatorial proof for this result.
 The tree record number are already present in the literature  in the context of queues
and traffic lights in the work of Haight  \cite{Haight}. 
 For more information on this sequence, visit \url{https://oeis.org/A259334}.

Observe that 
Corollary \ref{pablos_lemma}, implies that the record number $R_{\bullet}(n,k)$, counting the number of rooted trees of order $n$ and $k$ records, also counts  the number of connected endofunctions of girth $k$ minus the number of those of girth $k+1.$ 
The formula for tree record numbers is the obtained by combinating this observation with Cayley's forest formula \cite{Cayley}, an elegant 1889 result of Cayley stating that the number of rooted forests on $n$ vertices rooted at $[k]$ is equal to
$
     k n^{n-k-1}
$. We give a new proof for this elegant result in   Corollary  \ref{prop: bijection [k]-rooted} based on the notion of tree record.

\begin{cor}[Theorem 4.1 of \cite{LRT-GF-records}]
\label{prop:R_planted}
  The number of  rooted   trees labeled with  $[n]$ and  with $k$ records obeys the equation
    \[
    \planted(n,k) = 
       k \,(n-1)\cdots
    (n-k+1) \, n^{n-k-1}.
     \]
\end{cor}

\begin{proof}
The number of endofunctions of girth $k$ is 
\begin{align}
\label{endomophism_girth_k}
n(n-1)\dots(n-k+1) \ n^{n-k-1}
\end{align}
as $kn^{n-k-1}$ counts the number of  labeled rooted forests  with a fixed $k$-set of roots (Cayley's forest formula), and $n(n-1)\dots(n-k+1)/k$ counts the number of ways to choose the $k$-set of roots of the trees and  order it cyclically. Finally, Corollary \ref{pablos_lemma} allows us to deduce the formula for the record numbers from expression (\ref{endomophism_girth_k}) by subtracting the count of the number of endofunctions of girth $k+1$ from the count of those of girth $k$.
\end{proof}

\begin{cor}
Some subsequences of tree record numbers:
\begin{enumerate}
    \item The number of increasing trees on $n$ nodes is 
$\planted(n,n) = (n-1)!$.

\item The number of unrooted trees of order $n$ is $\planted(n,1) = n^{n-1}$.

\item The  number of  rooted trees in all rooted forests on n nodes, including the common root is $\planted(n+2,2) = 2(n+1)  (n+2)^{n-1}$. This is sequence \url{https://oeis.org/A089946}.

\item The sum of the positions of the records of all permutations of $\SS_n$ is 
$\planted(n+1,n) = n \cdot n!$. This is sequence \url{https://oeis.org/A001563}.

\end{enumerate}
\end{cor}
These statements can be directly checked by comparing the specific specializations of the record numbers with the formulas for the sequences appearing in the OEIS.
Other interesting subsequences can be found at  \url{https://oeis.org/A259334}.

\section{The genesis sequence of the online enclyclopedia of integer sequences}

\label{se:genesis}

In this section we the elegant connections between records in trees, connected endofunctions and the very first sequence in the OEIS database that we refered to in the introduction to this work.
The \emph{genesis sequence} $\Gamma = (\gamma_k)_{k}$ of the OEIS  is defined by setting $\gamma_k$ to be the sum of the height of all nodes of  all rooted trees of order $k$, divided by $k$. Due to this division, this sequence is referred to as the normalized total height of all nodes in all rooted trees with $n$ labeled nodes: \url{https://oeis.org/A000435}.

In this section, we show that the genesis sequence also counts non-root records in all rooted trees of order $n$.
Our main combinatorial tool for this section will be the notion of a catalyst.
Let $T$ be a labeled rooted tree. We define a \emph{catalyst} for $T$ to be a pair $(u, v)$ of distinct vertices of $T$ such that $v$ is an ancestor of $u$, and denote by $C(T)$ the set of catalysts for $T$.  A closely related, but different, notion of  Catalyst is introduced in \cite{briand2025nonintersectingpathsdeterminantdistance}.

\begin{lem}
\label{lem:catalyest_total_height}
Let $T$ be a rooted tree. The number of catalysts for $T$ is equal to the sum of the heights of all vertices in $T$. 
\end{lem}
We now state our main result. 
\begin{thm}
\label{thm:catalysts}
The  total number of non-root records in all rooted trees labeled with $[n]$ is equal to $\gamma_n$, the normalized total height of all nodes in all rooted trees with $n$ labeled nodes.
\end{thm}

Again, we present a bijective proof of this result. 
\begin{proof}
    Let $P(n)$ denote the set of pairs $(T, r)$ where $T$ is a rooted tree labeled with $[n]$ and $r$ is a record of $T$ different from the root. 
    We exhibit an explicit bijection between the set $[n]\times P(n)$ and the set \[
    C(n) = \bigsqcup_T C(T).
    \] 
    Note that the existence of such a bijection implies the result by Lemma \ref{lem:catalyest_total_height}. 
    Let $(v, T, r)$ be an element of $[n] \times P(n)$.
     We shall construct a tree $T'$ and a vertex $x$ such that $(x, v)$ is a catalyst for $T'$. 
     
    First, we compute the endofunction $f$ corresponding with the pair $(T, r)$ under the bijection of Theorem \ref{thm:connected endofunctions}. Let $\sigma$ denote the cycle of $f$, and
    let $x$ be the first element in the orbit of $v$ under $f$ that belongs to the cycle $\sigma$. Define an endofunction $g$ by
    \[
    g(z) = f(z) \quad \text{if} \ z\neq \sigma^{-1}(x), \qquad g(\sigma^{-1}(x)) = v.
    \]
    Let $\pi$ denote the cycle of $g$ and let $w = \pi^{-1}(x)$. 
    Lastly, define $T'$ to be the tree with parent map $p$, where $p$ is given by 
    \[
    p(z) = g(z) \quad \text{if } z\neq w, \qquad p(w)=w.  
    \]
    Then root $T'$  at $w$. It is simple to check that $(x, v) \in C(T')$, and that all steps are reversible.    
\end{proof}

\begin{figure}
    \centering
    $\centering
\vcenter{\hbox{
\begin{tikzpicture}[
    circ/.style={circle, draw, semithick, minimum size=1.5mm, inner sep=0pt, fill=black},
    arr/.style={-stealth, semithick}
]
    \def\sepFirstLevel{0.8}
    \def\len{0.8}
    \node[circ, label={above:3}] (3) at (0, 0) {};
    \node[circ, label={right:$4$}] (4) at ($(3) + (-50:\len)$) {};
    \node[circ, label={right:$5$}] (5) at ($(3) + (-130:\len)$) {};
    \node[circ, label={below:$6$}] (6) at ($(5) + (-120:0.9*\len)$) {};
    \node[circ, label={below:$2$}] (2) at ($(5) + (-60:0.9*\len)$) {};
    
    \node[circ, label={below:$1$}] (1) at ($(4) + (-120:0.9*\len)$) {};
    \node[circ, label={below:$7$}] (7) at ($(4) + (-60:0.9*\len)$) {};

    \draw (6) -- (5);
    \draw (2) -- (5);
    \draw (5) -- (3);
    \draw (4) -- (3);
    \draw (1) -- (4);
    \draw (7) -- (4);
    
    \draw[draw=USred, semithick, fill=white] (6) circle [radius=1.6mm];
    \node[circ, fill=USred, draw=USred,
    label={left:$\textcolor{USred}{v}$}
    ] (6bis) at (6) {};
    
    \draw[draw=darkgreen, semithick, fill=white] (4) circle [radius=1.6mm];
    \node[circ, fill=darkgreen, draw=darkgreen, label={left:$\textcolor{darkgreen}{r}$}] (4bis) at (4) {};
\end{tikzpicture}}}
\qquad
\vcenter{\hbox{
\begin{tikzpicture}[
    circ/.style={circle, draw, semithick, minimum size=1.5mm, inner sep=0pt, fill=black},
    arr/.style={-stealth, semithick}
]
    \def\sepFirstLevel{0.8}
\def\len{0.8}
    \node[circ, fill=darkblue, draw=darkblue, label={above:3}, label={below:$\textcolor{darkblue}{x}$}] (3) at (0, 0) {};
    \node[circ, draw=darkgreen, fill=darkgreen, label={above:$4$}, label={below:$\textcolor{darkgreen}{r}$}] (4) at ($(3) + (0:\len)$) {};
    \node[circ, label={above:$5$}] (5) at ($(3) + (-180:\len)$) {};
    \node[circ, label={above:$7$}] (7) at ($(4) + (30:\len)$) {};
    \node[circ, label={above:$1$}] (1) at ($(4) + (-30:\len)$) {};
    \node[circ, label={above:$6$}, label={left:$\textcolor{USred}{v}$}] (6) at ($(5) + (150:\len)$) {};
    \node[circ, label={above:$2$}] (2) at ($(5) + (210:\len)$) {};
    %
    

    \begin{scope}[-stealth]
    \draw (6) -- (5);
    \draw (2) -- (5);
    \draw (5) -- (3);
    \draw (4) edge[bend right] (3);
    \draw (3) edge[bend right] (4);
    \draw (1) -- (4);
    \draw (7) -- (4); 
    \end{scope}
    \draw[draw=USred, semithick, fill=white] (6) circle [radius=1.6mm];
    \node[circ, fill=USred, draw=USred] (6bis) at (6) {};
\end{tikzpicture}
}}
\qquad
\vcenter{\hbox{
\begin{tikzpicture}[
    circ/.style={circle, draw, semithick, minimum size=1.5mm, inner sep=0pt, fill=black},
    arr/.style={-stealth}
]
    \newcommand{\len}{0.8}
    \begin{scope}    
    \node[circ,  label={right:$\scriptsize 6$}, label={above:$\textcolor{USred}{v}$}] (c1) at (0, 0) {};
    \node[circ, draw=darkorange, fill=darkorange, label={right:$\scriptsize 5$}, label={below:$\textcolor{darkorange}{w}$}] (c2) at ($(c1) + (-135:\len)$) {};
    \node[circ,  label={right:$\scriptsize 3$}] (c3) at ($(c2) + (-45:\len)$) {};
    \node[circ, label={right:$\scriptsize 4$}] (c4) at ($(c3) + (45:\len)$) {};
    \end{scope}

    \node[circ, label={right:$\scriptsize 7$},] (15) at ($(c4) + (45:\len)$) {};
    \node[circ, label={above right:$\scriptsize 1$},] (11) at ($(c4) + (-45:\len)$) {};
    
    \node[circ, label={above:$\scriptsize 2$},] (10) at ($(c2) + (180:\len)$) {};

    \begin{scope}
    \draw[arr] (c1) -- (c2);
    \draw[arr] (c3) -- (c4);
    \end{scope}
    \draw[arr] (c4) -- ($(c1)+(-45:1.6mm)$);
    \draw[arr] (c2) -- ($(c3)+(135:1.6mm)$);

    \draw[arr] (10) -- (c2);
    \draw[arr] (11) -- (c4);
    \draw[arr] (15) -- (c4);
    
    \draw[draw=USred, semithick, fill=white] (c1) circle [radius=1.6mm];
    \node[circ, fill=USred, draw=USred] (c1bis) at (c1) {};
    
    \draw[draw=darkblue, semithick, fill=white] (c3) circle [radius=1.6mm];
    \node[circ, fill=darkblue, draw=darkblue, label={below:$\textcolor{darkblue}{x}$}] (c1bis) at (c3) {};
    
\end{tikzpicture}
}}
\qquad
\vcenter{\hbox{
\begin{tikzpicture}[
    circ/.style={circle, draw,  minimum size=1.5mm, inner sep=0pt, fill=black},
    arr/.style={-stealth}
]
    \def\len{0.8}
    \node[circ, draw=darkorange, fill=darkorange, label={right:$\scriptsize 5$}, label={above:$\textcolor{darkorange}{w}$}] (5) at (0, 0) {};
    \node[circ, label={right:$\scriptsize 2$},] (2) at ($(5) + (-60:\len)$) {};
    \node[circ, label={right:$\scriptsize 6$},] (6) at ($(5) + (-120:\len)$) {};
    \node[circ, label={right:$\scriptsize 4$},] (4) at ($(6) + (-90:\len)$) {};
    
    \node[circ, label={below:$\scriptsize 1$},] (1) at ($(4) + (-130:\len)$) {};
    \node[circ, label={below:$\scriptsize 3$},] (3) at ($(4) + (-90:\len)$) {};
    \node[circ, label={below:$\scriptsize 7$},] (7) at ($(4) + (-50:\len)$) {};

    \draw (1) -- (4);
    \draw (3) -- (4);
    \draw (7) -- (4);
    \draw (4) -- (6);
    \draw (6) -- (5);
    \draw (2) -- (5);
    
    \draw[draw=darkblue, semithick, fill=white] (3) circle [radius=1.6mm];
    \node[circ, fill=darkblue, draw=darkblue, label={left:$\textcolor{darkblue}{x}$}] (3bis) at (3) {};
    
    \draw[draw=USred, semithick, fill=white] (6) circle [radius=1.6mm];
    \node[circ, fill=USred, draw=USred,label={left:$\textcolor{USred}{v}$}] (6bis) at (6) {};
\end{tikzpicture}}}$
    \caption{The bijection of Theorem \ref{thm:catalysts}.}
    \label{fig:first_sequence_bijection}
\end{figure}
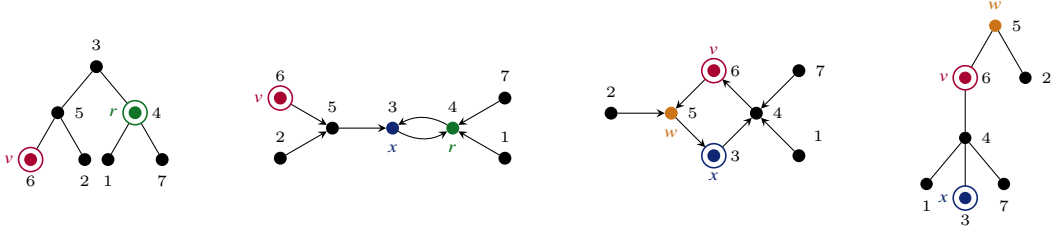

An explicit example of the bijection of Theorem \ref{thm:catalysts} is illustrated in Figure \ref{fig:first_sequence_bijection}. 

\begin{cor} The genesis sequence can be defined in terms of the record numbers setting
    \[
    \gamma_k =
    \sum_{k=1}^n (k-1)R_\bullet(n, k).
    \]
\end{cor}

As recounted by Sloan in his podcast, Riordan and Sloan \cite{Riordan_Sloane}, 
obtained an elegant expression for the normalized height of all nodes in all trees of order $n$. Let $W(z)$ be the generating function for trees according to their total height (without any normalization). The key element of
Riordan and Sloan's formula for the coefficients of $W(z)$ is a beautiful   manipulation of generating functions from which they derive the functional equation $
W(z) = \left( z \mathcal{T}'(z)\right)^2$. For the sake of completeness, we give a simple bijective proof of the symbolic version of this identity.

\begin{lem}[Riordan-Sloan, \cite{Riordan_Sloane}]
\label{lem:Riordan_Sloane}
There exists a bijection between the set of rooted trees and a catalyst defined on them, and the set of pairs of rooted tree with a selected node on each.
\end{lem}

For the proof of this result we present an explicit bijection.  

\begin{proof}
By Lemma \ref{lem:catalyest_total_height}, $W(z)$ is the generating function for rooted trees with a catalyst. On the other hand, $\left( z \mathcal{T}'(z)\right)^2$  is the generating function for pairs of rooted trees with a selected node in each of them. 

Starting from a pair of rooted trees with a selected node in each of them, we construct a rooted tree $T$ and a catalyst for $T$. To do so,  we add an edge between the two roots, designating the selected node of the first tree as the root of $T$, and define the catalyst  to be the pair $(a,b)$ where $a$ is the root of the first tree and  $b$ is the distinguished node of the second one. 
To see that this is a bijection, we observe that the root $r$ of the second tree can be characterised as the  child of $a$ that belongs to the unique path connecting $a$ to $b.$ Thus, we recover the initial pair by deleting the edge $(a,r)$. After this, it is immediate to check that this process can be reversed.

\end{proof}

\begin{figure}
    \[\begin{tikzpicture}[
    circ/.style={circle, draw, semithick, minimum size=1.5mm, inner sep=0pt, fill=black},
    hl/.style={circ, fill=USred, draw=USred},
    hl2/.style={circ, fill=darkblue, draw=darkblue}
]
    \newcommand{\len}{0.8}
    \node[hl2, label={left:$5$}] (root1) at (0, 0) {};
    \node[circ, label={left:$3$}] (A) at ($(root1) + (-120:\len)$) {};
    \node[circ, label={right:$2$}] (B) at ($(root1) + (-60:\len)$) {};
    \node[circ, label={left:$11$}] (F) at ($(A) + (-90:\len)$) {};
    \node[circ, label={left:$7$}] (C) at ($(F) + (-120:\len)$) {};
    \node[circ, label={right:$1$}] (D) at ($(F) + (-60:\len)$) {};
    \node[circ, label={left:$10$}] (E) at ($(C) + (-90:\len)$) {};

    \draw (root1) -- (A);
    \draw (root1) -- (B);
    \draw (F) -- (A);
    \draw (C) -- (F);
    \draw (D) -- (F);
    \draw (E) -- (C);
    
    \node[hl2, label={right:$8$}] (root2) at (2, 0) {};
    \node[circ, label={left:$12$}] (G) at ($(root2) + (-120:\len)$) {};
    \node[circ, label={right:$9$}] (K) at ($(root2) + (-60:\len)$) {};
    \node[circ, label={below:$6$}] (H) at ($(G) + (-130:\len)$) {};
    \node[circ, label={below:$13$}] (I) at ($(G) + (-90:\len)$) {};
    \node[circ, label={below:$4$}] (J) at ($(G) + (-50:\len)$) {};

    \draw (H) -- (G);
    \draw (I) -- (G);
    \draw (J) -- (G);
    \draw (G) -- (root2);
    \draw (K) -- (root2);

    \draw[semithick, dashed, darkblue] (root1) -- (root2);

    \draw[draw=USred, semithick, fill=white] (I) circle [radius=1.6mm];
    \draw[draw=USred, semithick, fill=white] (A) circle [radius=1.6mm];
    \node[hl, USred] (Ibis) at (I) {};
    \node[hl, USred] (Abis) at (A) {};
\end{tikzpicture}
    \qquad
    \begin{tikzpicture}[
    circ/.style={circle, draw, semithick, minimum size=1.5mm, inner sep=0pt, fill=black},
    hl/.style={circ, fill=USred, draw=USred},
    hl2/.style={circ, fill=darkblue, draw=darkblue}
    ]
    \def\sepFirstLevel{1}
    \def\len{0.8}
    \node[hl, label={right:$3$}] (3) at (0, 0) {};
    \node[hl, label={right:$5$}] (4) at ($(3) + (-45:1.1*\len)$) {};
    \node[circ, label={right:$11$}] (5) at ($(3) + (-135:1.1*\len)$) {};
    \node[circ, label={left:$7$}] (6) at ($(5) + (-120:\len)$) {};
    \node[circ, label={left:$1$}] (2) at ($(5) + (-60:\len)$) {};
    
    \node[circ, label={right:$2$}] (1) at ($(4) + (-120:\len)$) {};
    \node[hl2, label={right:$8$}] (7) at ($(4) + (-60:\len)$) {};

    \node[circ, label={left:$10$}] (E) at ($(6) + (-90:\len)$) {};

    \node[circ, label={right:$9$}] (K) at ($(7) + (-60:\len)$){};
    \node[circ, label={right:$12$}] (G) at ($(7) + (-120:\len)$){};

    \node[circ, label={below:$6$}] (H) at ($(G) + (-130:\len)$) {};
    \node[circ, label={below:$13$}] (I) at ($(G) + (-90:\len)$) {};
    \node[circ, label={below:$4$}] (J) at ($(G) + (-50:\len)$) {};
    
    \draw (6) -- (5);
    \draw (2) -- (5);
    \draw (5) -- (3);
    \draw (4) -- (3);
    \draw (1) -- (4);
    \draw[semithick, darkblue] (7) -- (4);
    \draw (G) -- (7);
    \draw (K) -- (7);
    \draw (H) -- (G);
    \draw (I) -- (G);
    \draw (J) -- (G);
    \draw (E) -- (6);
    
    \draw[draw=USred, semithick, fill=white] (I) circle [radius=1.6mm];
    \draw[draw=darkblue, semithick, fill=white] (4) circle [radius=1.6mm];
    \node[hl, USred] (Ibis) at (I) {};
    \node[hl, darkblue] (4bis) at (4) {};
\end{tikzpicture}
    \]
    \caption{On the left, two rooted trees with a distinguished vertex each. On the right, the corresponding catalyst under the bijection of Lemma \ref{lem:Riordan_Sloane}.}
    \label{fig:Riordan_Sloane_bijection}
\end{figure}
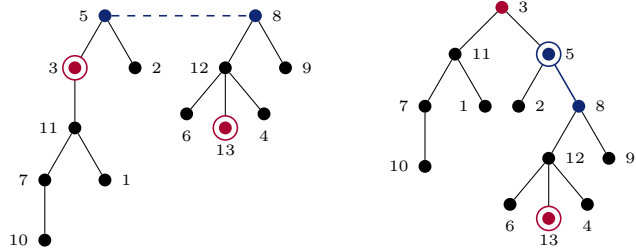

Finally, from the functional equation $
W(z) = \left( z \mathcal{T}'(z)\right)^2$ (Lemma \ref{lem:Riordan_Sloane}) together with Abel's
generalization of the binomial formula, Riordan and Sloan obtain the expression for $\gamma_i$ at the right-hand-side of the equation of Eq.~(\ref{eq:R-S}).
\begin{cor}[Riordan-Sloan and the record numbers] 
The number of connected endofunctions without a fixed point is 
  \begin{align}
  \label{eq:R-S}
  \gamma_k = \sum_{k=1}^n (k-1) \planted(n,k) = (n-1)! \sum_{k=0}^{n-1} \frac{n^k}{k!}
    \end{align}
\end{cor}
Note that combining Eq.~(\ref{eq:R-S}) with the record number formula, we obtain that both sides of the following nontrivial identity count the number of connected endofunctions.
\begin{cor}[Connected endofunctions and the record numbers]
\label{Cor:records_connected_endofunctions}
The number of connected endofunctions from $[n]$ to $[n]$ is given by both sides of the following identity
   \begin{align}  \label{Eq:records_connected_endofunctions} 
 \sum_{k=1}^n
       k^2 \, \frac{(n-1)!}{(n-k)!} \, n^{n-k-1}
       \, = \,
            n^{n-1} + (n-1)! \sum_{j=0}^{n-2} \frac{n^j}{j!}.
   \end{align}
   \end{cor}
 Note that expression for the number of connected endofunctions on the right-hand side of Eq.~(\ref{Eq:records_connected_endofunctions}), the one that doesn't involve records, appears in the work of \cite{Katz} on his study of probabilistic aspects of connected endofunctions. Moreover, its first summand counts connected endofunctions with a fixed point (that we identify with the root). Thus, $(n-1)! \sum_{j=0}^{n-2} \frac{n^j}{j!}$ counts the number of connected endofunctions with no fixed points.\\


\subsection{The weight of the symmetric group}

The \emph{weight of $\SS_n$}, denoted by $W(n)$, is the sum of the reflection lengths of all permutations of $n$ letters. where the reflection length of a permutation $\pi$ is defined as the minimal number of transpositions whose product equals $\pi$.  This is sequence \url{https://oeis.org/A067318} of the OEIS. For example, the weight of $\SS_2$ is seven, as $\emptyset$, $(12)$, $(13)$, $(23),$ $(12)(23)$ and $(13)(23)$ is a list of all of its elements written in a minimal way as a product of transpositions.
 
This subsection is devoted to the proof of the following result. 

\begin{p}
\label{prop:weightSn}
The forest record number $R(n,n-1)$ is equal to the weight of $\SS_n$ for every  natural number $n$.
\end{p}

Let $\mathcal R_m(n,k)$ be the set of rooted forests labelled with $[n]$, with $m$ connected components and with $k$ records, and let $\mathcal R(n, k) = \bigsqcup_m \mathcal R_m(n,k)$.
There is a well-known bijection \cite{Stanley_EC1} between the set $\mathcal R_m(n, n)$ of (rooted) increasing forests with $m$ trees, and the set of permutations of $\mathbb S_n$ with $m$ left-to-right minima. The cardinality of both families is given by the \emph{unsigned Stirling numbers of the first kind, $c(n, m)$}, which also count the permutations of $\mathbb S_n$ with $m$ cycles. The latter claim is a direct consequence of Foata's fundamental transformation \cite{Foata}. 

A minimal factorization of a permutation in $\mathbb S_n$ with $m$ cycles consists of $n - m$ transpositions. On the other hand,  there are $c(n, m)$ such permutations. This allows us to express the weight of the symmetric group 
in terms of the Stirling numbers: 
\begin{equation}
\label{eqn: weight Sn stirling}
W(n) = \sum_{m=1}^n (n-m)c(n, m).
\end{equation}
Notice that summand $(n - m)c(n, m)$ also gives the number of non-root vertices in all increasing forests with $m$ trees and $n$ vertices. Indeed, there are $c(n, m)$ such forests, and each one has $n - m$ non-root vertices. 
Let $\mathcal N_m(n)$ be the set of pairs $(v, F)$ where $F$ is a forest in $\mathcal R_m(n, n)$ and $v$ is a non-root vertex of $F$. 
Equation \eqref{eqn: weight Sn stirling} can be rewritten as 
\[
W(n) = 
\sum_{m=1}^n |\mathcal N_m(n)|.
\]
To complete the proof of Proposition \ref{prop:weightSn}, we present, in Lemma \ref{lem:equinumerous}, a bijection  between $\mathcal N_m(n)$ and $\mathcal R_m(n, n-1)$. This allows us to conclude 
\[
W(n) = \sum_{m=1}^n |\mathcal R_m(n, n-1)| = |\mathcal R(n,n-1)| = R(n, n-1),
\]
because the sets $\mathcal R_m(n, n-1)$ partition $\mathcal R(n, n-1)$.

\begin{figure}[h!]
\centering
\scalebox{0.9}{
\begin{tikzpicture}[
    level distance=0.7cm,
    level 1/.style={sibling distance=0.9cm},
    level 2/.style={sibling distance=0.75cm},
    level 3/.style={sibling distance=0.75cm},
    circ/.style={circle, draw, semithick, minimum size=1.5mm, inner sep=0pt, fill=black},
    hl/.style={circ, fill=USred, draw=USred},
]\matrix(outer) [matrix of nodes, row sep=0.3cm, column sep=1cm, ampersand replacement=\&]{ 
    \node[record, label={right:$\circ$}] (N-0) {}child {	
	    node[record, label={right:1}] (N-1) {}child {	
		    node[record, label={right:2}] (N-2) {}child {	
			    node[record, label={right:5}] (N-5) {}}}child {	
		    node[record, label={right:3}] (N-3) {}}child {	
		    node[record, label={right:4}] (N-4) {}child {	
			    node[record, label={right:6}] (N-6) {}child {	
				    node[record, label={right:7}] (N-7) {}}}child {	
			    node[record, label={right:8}] (N-8) {}}child {	
			    node[record, label={right:9}] (N-9) {}}}}; 

    \phantom{\node[circle,draw=USred, semithick, fill=white, label={below:$\textcolor{USred}{v}$}] (c) at (N-8){};}; \&
    \draw[->] (-1, -2) -- (0,-2) node [midway, above] {$\phi$}; \&
    \node[record, label={right:$\circ$}] (N-0) {}child {	
	    node[record, label={right:1}] (N-1) {}child {	
		    node[record, label={right:2}] (N-2) {}child {	
			    node[record, label={right:5}] (N-5) {}}}child {	
		    node[record, label={right:3}] (N-3) {}}child {	
		    node[record, label={right:6}] (N-6) {}child {	
			    node[record, label={right:7}] (N-7) {}}child {	
			    node[record, label={right:8}] (N-8) {}child {	
				    node[non-record, label={right:4}] (N-4) {}child {	
					    node[record, label={right:9}] (N-9) {}	
						     }}}}};\\};
    
    \draw[draw=USred, semithick, fill=white] (c) circle [radius=1.6mm];
    \node[hl, USred, label={below:$\textcolor{USred}{v}$}] (Ibis) at (c) {};
\end{tikzpicture}}

\caption{The bijection $\phi$. In this example, $v = 8$, $l = 1$ and $(x_0, x_1, x_2, x_3) = (g, x_1, v, p) = (1,6, 8, 4)$.
}
\label{fig:weight_Sn_bijection}
\end{figure}
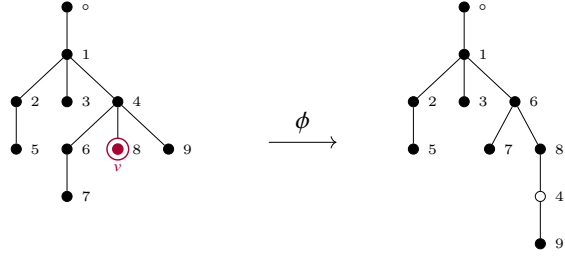

\begin{lem}
\label{lem:equinumerous}
    There exists a bijection 
    $\phi: \mathcal N_m(n) \to \mathcal R_m(n, n-1)$. 
\end{lem}
\begin{proof}
Take a pair $(v, F)$ in $\mathcal N_m(n)$, and attach all roots of $F$ to a common parent $\circ$, and let $T$ be the resulting tree, rooted at $\circ$. 

Let $p$ and $g$ denote the parent and grandparent of $v$ in $T$, and let $x_1 < x_2 < \cdots< x_l$ be the children of $p$ that are smaller than $v$. 
For convenience, we also write $x_0 = g$,  $x_{l + 1} = v$ and $x_{l + 2} = p$. We construct a tree $T'$ from $T$ as follows: 
\begin{enumerate}
    \item[--] The parent of $x_i$ in $T'$ is $x_{i-1}$, for all $i \in [l+2]$.
    \item[--]  The parent of any other vertex is the same as in $T$. 
\end{enumerate}

Finally, we define $F' =\phi(v, T)$ as the rooted forest obtained by removing the virtual vertex $\circ$ from $T'$. See Figure \ref{fig:weight_Sn_bijection}. It is easy to verify that all vertices in $F'$ are records, with the sole exception of $p$.

Let us show that $\phi$ is a bijection. Let $F' \in \mathcal R_{m}(n, n-1)$ and let $p$ be the unique vertex of $F'$ that is not a record. Once the vertices $x_0, x_1, \dots, x_{l + 2}$ have been determined, finding the inverse map is straightforward. This turns out to be quite simple. First, we add back the virtual vertex $\circ$. Then, we identify the first vertex smaller than $p$ in the path from $p$ towards the root; this vertex is $x_{0}$. All other vertices $x_i$ appear in order in the path from $x_0$ towards $x_{l+2} = p$.\qedhere
\end{proof}
\begin{cor}
    For $n \geq 1$ we have the equality 
    \[
     R(n, n-1) = \sum_{m= 1}^n (n-m)c(n, m)
    \]
where $c(n,m)$ denotes unsigned Stirling numbers of the first kind.
\end{cor}

\section{On the record  code  as an endofunction.} 
\label{se:record_code_endofunction}
\subsection{The record  code of a rooted tree.} 

To introduce the record code of an endofunction, we first state some basic ideas of \cite{LRT-Weary}.
A \emph{bonsai tree} is a labelled rooted tree whose only record is its root. Bonsai trees are labelled by natural numbers.
The \emph{bonsai sequence} of $T$ is the set of maximal bonsai trees that are subgraphs of $T$, ordered from smallest to largest by their record labels. 
These bonsais are obtained by deleting the edges connecting each of the records and their respective parents. The corresponding sequence of parents of the records, ordered also from smallest to largest, and starting at the second-smallest record is the \emph{attachment sequence}. We omit the parent of the smallest record in this sequence, since it must be $\circ$. The pair formed by the bonsai sequence  and the attachment sequence defines the \emph{record decomposition} of $T$. The record decomposition of $T$ is an alternative way of representing the tree $T$, since one may reconstruct $T$ from it.

The record decomposition may be expressed in a compact way akin to the better-known Pr\"ufer code. Recall that we denote  by $p_T$ the \emph{parent function} of $T$, which sends each node to its parent, and $\circ$ to itself. The \emph{record code} of $T$ is the sequence $r = (r_1,\dots,r_{n-1}) \in \{0,\dots,n\}^{n-1}$ for which $r_i = p_T(i)$ if $i$ is a non-record of $T$, and $r_i = p_T(i')$ if $i$ is a record of $T$, where $i'$ is the next largest record of $T$.

The record code is due to \cite{Kreweras-Moszkowski, Picciotto}. Each tree has a unique record code, and each sequence $r \in \{0,\dots,n\}^{n-1}$ corresponds to a unique tree. The record code therefore gives a bijective proof of Cayley's rooted formula. 

Finally, we add a note on the decoding of the record code, which is of practical use in understanding the record code, and also something that we exploit in Section \ref{sec:doublycovering-records}. It is clear that if one wishes to decode the record code that one first needs to identify which elements of $[n-1]$ will be records of the associated tree $T$. For each $i \in [n-1]$, define $r(i) = r_i$, and set $r(0) = 0, r(n) = n$. It is not difficult to see that the records of $T$ are exactly those elements $i \in [n]$ which are at least as large as $r^m(i)$ for each positive integer $m$. If $i$ is not a record, then since there is a record on the path to the root, some $r^k(i)$ must be larger than $i$. On the other hand, if $i$ is a record, then $r(i)$ is the parent of the next largest record $i'$. It is not difficult to show that any element of the form $r^k(i)$ must then either be a record $\leq i$ or a strict ancestor (i.e., an ancestor of the node that is not the node itself) of one of the records $\leq i'$. Each such element can be at most $i$.

\begin{ex} \label{ex:record_code}
A rooted tree $T$ (with record nodes drawn in black) and its record code:
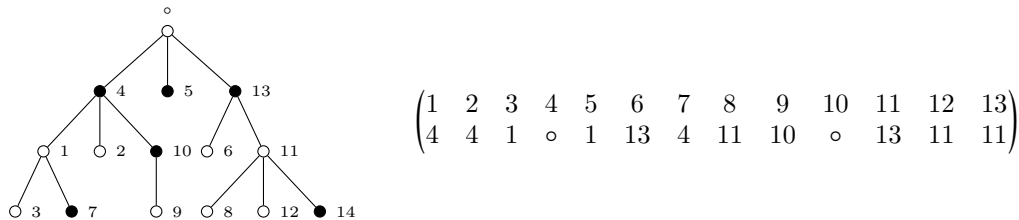
\begin{figure}[h!]
\centering
\begin{tabular}{@{} c c @{}}
\normalfont
\begin{tikzpicture}[
    level distance=0.8cm,
    level 1/.style={sibling distance=0.9cm},
    level 2/.style={sibling distance=0.75cm,},
    level 3/.style={sibling distance=0.75cm},
    ]
    \node[non-record, label={$\circ$}] (1) {}
    child {
    	node[record, label={right:4}] (5) {}
    	child {
    		node[non-record, label={right:1}] (2) {}
    		child {
    			node[non-record, label={right:3}] (4) {}
    		}
    		child {
    			node[record, label={right:7}] (8) {}
    		}
    	}
    	child {
    		node[non-record, label={right:2}] (3) {}
    	}
    	child {
    		node[record, label={right:10}] (11) {}
    		child {
    			node[non-record, label={right:9}] (10) {}
    		}
    	}
    }
    child {
    	node[record, label={right:5}] (6) {}
    }
    child {
    	node[record, label={right:13}] (14) {}
    	child {
    		node[non-record, label={right:6}] (7) {}
    	}
    	child {
    		node[non-record, label={right:11}] (12) {}
    		child {
    			node[non-record, label={right:8}] (9) {}
    		}
    		child {
    			node[non-record, label={right:12}] (13) {}
    		}
    		child {
    			node[record, label={right:14}] (15) {}
    		}
    	}
    };

\node (recordcode) at (7.3, -1.2) {\small $\begin{pmatrix}
1 & 2 & 3 & 4 & 5 & 6 & 7 & 8 & 9 & 10 & 11 & 12 & 13\\
4 & 4 & 1 & \circ & 1 & 13 & 4 & 11 & 10 & \circ & 13 & 11 & 11
\end{pmatrix}$};
\end{tikzpicture}
&
\end{tabular}
\caption{A rooted tree and its record code.}
\label{fig:example_records}
\end{figure}

To decode this record code, we first identify the records. For example, we check that $7$ is a record, since $r(7) = 4, r^2(7) = \circ$, and $r^i(7) = \circ$ for each $i \geq 3$, and so $7 \geq r^k(7)$ for each positive integer $k$. Similarly, we may check that $12$ is a non-record since $r^2(12) = 13$ is larger than $12$.

Once we have identified the records, we define the parent function $p_T$ as follows. Set $p_T(i) = r(i)$ for each non-record $i$, set $p_T(j) = \circ$ for the smallest record $j$, and set $p_T(k) = r(k')$ for each other record $k$, where $k'$ is the next smallest record to $k$.

\end{ex}

\subsection{On the record code and endofunctions} \label{sec:doublycovering-records}

Let  $T$ be a rooted tree $[n]_0$. We consider the record code of $T$ as a function from $[n-1]$ to $[n]$, and extend it into a  endofunction $\varphi_T :[n]_0 \to [n]_0$ by associating to each element $i$ in $[n-1]$ its image according to the  record 
code, (that is the $i^{th}$ entry of the record code), and setting $\varphi_T(\circ) = \circ,$ and $ \varphi_T(n) = n$. We are interested in $\varphi_T$ as it is exactly the endofunction we defined for the decoding process of the record code. 

Since the fixed points $(\circ)$ and $(n)$ are always cycles of $\varphi_T$, we distinguish them from the other cycles of $\varphi_T$, by referring to $(\circ)$ and $(n)$ as the \emph{trivial cycles} of $\varphi_T$, and to all other cycles as the \emph{non-trivial cycles} of $\varphi_T$.

Endofunctions on $[n]_0$ with fixed points $\circ$ and $n$ (that is, containing both trivial cycles) are in bijection with rooted trees $T$ on $[n]_0$.  To match with our definition of records of rooted trees on $[n]_0$, let us insist that $\circ$ is not a record of the endofunction. In such a way, the records of the endofunction $\varphi_T$ are exactly the records of $T$.
 
\begin{lem}
    The number of rooted trees on $[n]_0$ with $k$ records is equal to the number of endofunctions on $[n]_0$ with $\circ$ and $n$ as fixed points, and with $k$ records.
\end{lem}

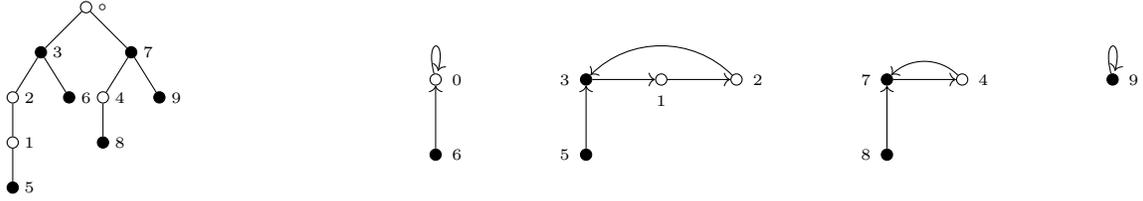
\begin{figure}[h!]
    \begin{minipage}{0.1\textwidth}
    \centering
    \begin{tikzpicture}
    [
    level distance=0.6cm,
    level 1/.style={sibling distance=1.2cm},
    level 2/.style={sibling distance=0.75cm,},
    level 3/.style={sibling distance=0.75cm},
    level 4/.style={sibling distance=0.75cm},
    cutedge/.style={dashed},
    every label/.style={inner sep=2pt,font=\tiny}
    ]
    \node[non-record, label={right:$\circ$}] (V-0) {}child {	
	    node[record, label={right:3}] (V-3) {}child {	
		    node[non-record, label={right:2}] (V-2) {}child {	
			    node[non-record, label={right:1}] (V-1) {}child {
                    node[record, label={right:5}] (V-5) {}
                }
                }}child {	
		    node[record, label={right:6}] (V-6) {}}}child {	
	    node[record, label={right:7}] (V-7) {}child {	
		    node[non-record, label={right:4}] (V-4) {}child {	
			    node[record, label={right:8}] (V-8) {}}}child {	
		    node[record, label={right:9}] (V-9) {}}};
\end{tikzpicture}
    \end{minipage}
    \begin{minipage}{0.25\textwidth}
    \
    \end{minipage}
    \begin{minipage}{0.4\textwidth}
        \centering
        \usetikzlibrary{automata}
\begin{tikzpicture}[circ/.style={circle, draw, semithick, minimum size=1.5mm, inner sep=0pt, fill=black},
    arr/.style={-stealth, semithick}
]

\node[non-record][label={right:$0$}] (0) at (0,0) {};
\node[record][label={right:$6$}] (6) at (0,-1) {};

\node[record][label={left:$3$}] (3) at (2,0) {};
\node[non-record][label={below:$1$}] (1) at (3,0) {};
\node[non-record][label={right:$2$}] (2) at (4,0) {};
\node[record][label={left:$5$}] (5) at (2,-1) {};

\node[record][label={left:$7$}] (7) at (6,0) {};
\node[non-record][label={right:$4$}] (4) at (7,0) {};
\node[record][label={left:$8$}] (8) at (6,-1) {};

\node[record][label={right:$9$}] (9) at (9,0) {};

\path (0) edge [loop above] node {} (0);
\path (9) edge [loop above] node {} (9);

\draw[->] (6)--(0);
\draw[->] (3)--(1);
\draw[->] (1)--(2);
\draw[->] (7)--(4);
\draw[->] (8)--(7);
\draw[->] (2) to[bend right=45] (3);
\draw[->] (5)--(3);
\draw[->] (4) to[bend right=45] (7);
\end{tikzpicture}
    \end{minipage}
    \hfill
    \caption{A tree (left) $T$ on $[9]_0$. Its records are $3,5,6,7,8,9$, and its doubly record-covering pairs are $(3,5)$ with distance $3$, and $(7,8)$ with distance $2$, and so $\delta(T) = (2,3)$. The record code of $T$ is $23173\circ47$. The endofunction (right) $\varphi_T : [9]_0 \to [9]_0$ corresponding to $T$. The two non-trivial cycles (corresponding to the two consecutive record pairs of $T$) are $(3,1,2)$ and $(7,4)$, and the specification of $\varphi_T$, $\lambda(T) = (2,3)$, is equal to the distance partition $\delta(T)$.}
    \label{fig:doublycovering-records}
\end{figure}

We denote the set of endofunctions on $[n]_0$ with fixed points $\circ$ and $n$ by $\mathcal{E}_n$. The \emph{type} of an endofunction $\varphi$ is the partition whose parts are the lengths of the cycles of $\varphi$ (noting that a fixed point has length $1$). We will denote the type of the endofunction $\varphi$ by $\tau(\varphi)$. Let $\alpha = 1^{m_1}2^{m_2}\dots$ and $ \beta = 1^{m'_1}2^{m'_2}\dots$ be partitions such that $\alpha \subseteq \beta$. That is, such that $m_i \geq m'_i$ for each positive integer $i$. By $\alpha \setminus \beta$ we denote the partition $1^{m_1 - m'_1}2^{m_2 - m'_2}\dots$. Thus, for $\varphi \in \mathcal{E}_n$, $\tau(\varphi) \setminus 11$ indicates the partition recording the lengths of the non-trivial cycles of $\varphi$.

Define a \emph{partial order $<_T$} on $[n]_0$ by setting $x <_T y$ if and only if $x$ is an ancestor of $y$ in $T$. Let $R$ and $R'$ be records of $T$. Then, $R$ is said to be \emph{doubly record-covering}  $R'$, written \emph{$R \lhd R'$},  if the following two conditions are satisfied (i) $R <_T R'$, and (ii) $R < R'$ and there is no record $R''$ such that $R < R''< R'$.
In words, $R \lhd R'$ if $R'$ is the next largest record to $R$ and is also a descendant of $R$. We call $\{R,R'\}$ a \emph{doubly record-covering pair}. 

We have elected to use the phrase ``doubly record-covering'' since one can characterize the two conditions via the condition ``$R$ covers $R'$'' in the following two posets both of whose ground-set is the set of records of $T$: (i)  the poset with the relation $<_T$, and (ii)  the poset with the relation $>$ on the labels.

Recall that for nodes $a,b$ in $T$, the distance $d_T(a,b)$ is the length of the unique path between $a$ and $b$.
 The \emph{distance partition $\delta(T)$} of $T$ is defined as  be the partition  whose parts are the distances $d_T(R,R')$ between doubly record-covering record pairs $R \lhd R'$  of $T.$ The weight  of $\delta(T)$ is 
 $d(T) = \sum_{R \lhd R'} d_T(R,R')$ and the  length of $\delta(T)$ is the number of doubly record-covering pairs of $T$.

In the next lemma we prove that the doubly record-covering pairs of $T$ correspond exactly to the non-trivial cycles of $\varphi_T$, such that the distances between the two records making up each doubly record-covering pair are exactly the lengths of the non-trivial cycles.

\begin{lem}
    Let $T$ be a rooted tree on $[n]_0$ with $\ell$ doubly record-covering pairs, and distance partition $\delta(T)$. 
    Then $\tau(\varphi_T) \setminus 11 = \delta(T)$.
\end{lem}

\begin{proof}

Consider a non-trivial cycle $C = (v_1,\dots, v_m)$ of the endofunction $\varphi_T$ so that $v_j = \varphi_{j+1}(v_1)$ for each $j=2,\dots,m$. We know that such a cycle can contain only a single record $R$ of the endofunction, and since the records of $\varphi_T$ are the same as that of $F$, $R$ is the only record of $F$ in $C$. Therefore, we may assume that $v_1 = R$ is the unique record of $C$. 

Let us now consider what the image of such a non-trivial cycle in $F$. Since $R$ is a record, we see that $v_2 = \varphi_T(v_1)$ is the parent of the record $R'$ in $F$ (where $R'$ is the next largest record to $R$). Moreover, since $R$ is the only record in $C$, we know that $\varphi_T(v_j)$ is simply the parent of $v_j$ in $F$ for each $j \neq 1$. Therefore, $R'$ must be a descendant of $R$, and since it is the next largest record to $R$, it follows that $R \lhd R'$. 
By the same token, whenever $R \lhd R'$, it is clear that we obtain a cycle $C$ of the endofunction whose length is $d_T(R, R')$. 
\end{proof}

An immediate consequence of this lemma is that the number of rooted trees on $[n]_0$ with $\ell$ doubly record-covering pairs  is equal to the number of endofunctions in $\mathcal{E}_n$ with $\ell$ non-trivial cycles. Similarly, the number of trees $T$ with distance partition $\delta(T)$ is equal to the number of endofunctions $\varphi \in \mathcal{E}_n$ satisfying $\tau(\varphi) \setminus 11 = \delta(T)$. 

\begin{ex}
    Consider the rooted tree $T$ in Figure \ref{fig:doublycovering-records} (left) with corresponding endofunction $\varphi_T$ (right).
    The two  doubly record-covering pairs $(3,5)$ and $(7,8)$ of $T$ correspond to the two non-trivial cycles of $\varphi_T$: $(3,1,2)$ and $(7,4)$. Then $\delta(T)$ is $(3,2)$, and the type of $\varphi_T$ is $(1,1,2,3)$, so $\delta(T) = \tau(\varphi_T) \setminus 11$. 
\end{ex}

We may now use the previous lemma in order to enumerate the number of trees on $[n]_0$ with $\ell$ doubly record-covering pairs, and also refine this count to track the distances between the doubly record-covering pairs.  

\begin{lem}
    Let $n$ and $q$ be positive integers with $q < n$, and let $\lambda = 1^{m_1}2^{m_2}\dots(n-1)^{m_{n-1}}$ be a partition of $q$ with $m_1 \geq 2$.
    The number of rooted trees $T$ on $[n]_0$ with distance partition $\delta(T) = 1^{m_1-2}2^{m_2}\dots (n-1)^{m_{n-1}} = \lambda \setminus 11$ is 
    \begin{equation}
    \frac{(n-1)!}{(n-|\lambda|+1)!}\
    \frac{|\lambda|(n+1)^{n-|\lambda|}}{z_{\lambda}}\frac{m_1!}{(m_1-2)!}.
    \end{equation}
    where $z_{\lambda} = \prod_{p\ge1} p^{m_p}\, m_p!$. 
\end{lem}

\begin{proof}

We enumerate the number of endofunctions $\varphi$ in $\mathcal{E}_n$ of type $\tau(\varphi) = \lambda$. Note that by definition $|\lambda| = q$ is the number of elements appearing in some cycle of $\varphi$. There are $\binom{n-1}{q-2}$ choices for non-trivial cycle vertices, and  
    $
    \frac{q!}{z_{\lambda}}\frac{m_1!}{(m_1-2)!}$
    ways to form the non-trivial cycles (these are simply the number of permutations of $[q]$ whose cycle type corresponds to the non-trivial cycles).   Finally, there are $q(n+1)^{n-q}$ choices for the remaining nodes (these are rooted forests on $\{0,1\dots,n\}$ with $q$ connected components). 
\end{proof}

\begin{cor}
    The number of rooted trees on $[n]_0$ with $\ell$ doubly record-covering pairs record pairs is 
\begin{equation} \label{eq:doublycovering}
    \sum_{q=\ell+2}^{n+1}\binom{n-1}{q-2}c(q-2,\ell)q(n+1)^{n-q}.
\end{equation}
\end{cor}

\begin{proof}
    If we are interested only in the number of cycles of the endofunction and not the type (and thus only the number of doubly record-covering pairs pairs, and not their corresponding distances), then the $\frac{(q-2)!}{z_{\lambda}}\frac{m_1!}{(m_1-2)!}$ term simply becomes simply the unsigned Stirling number of the 1st kind, $c(q-2,\ell)$,
    where $\ell$ is the total number of non-trivial cycles. Finally, we sum over all possible choices of $q$.
\end{proof}

The triangular array of numbers obtained from Eq.~\eqref{eq:doublycovering} is illustrated in Table \ref{tab:doubly-covering}. At this point, it does not seem to be in the \verb|oeis|. 

\begin{table}[h!] 
\centering
\renewcommand{\arraystretch}{1.3} 

    \begin{tabular}{c|cccccccccc}
        \toprule
$n \backslash \ell$ & $0$ & $1$ & $2$ & $3$ & $4$ & $5$ & $6$ 
        \\
        \midrule
$1$ & $1$   \\
$2$ & $2$   & $1$   \\
$3$ & $8$   & $7$   & $1$   \\
$4$ & $50$  & $59$  & $15$  & $1$  \\
$5$ & $432$ & $622$ & $215$ & $26$ & $1$ \\
$6$ & $4802$ & $8009$ & $3390$ & $565$ & $40$ & $1$ & \\
$7$ & $65536$ & $122696$ & $60214$ & $12415$ & $1225$ & $57$ & $1$ \\
        \bottomrule
    \end{tabular}
        \caption{The number of rooted trees on $[n]_0$ with $\ell$ doubly record-covering pairs.} 
        \label{tab:doubly-covering}
\end{table}

Counting rooted trees on $[n]_0$ with no doubly record-covering pairs record pairs, the formula \eqref{eq:doublycovering} simplifies to 
\[
2(n+1)^{n-2}
\]
which is the sequence https://oeis.org/A007334. This value counts facets in permutahedra \cite{Bostan}, ``metered parking functions'' \cite{Daugherty}, as well as other other combinatorial objects (see the references in the oeis entry). However, perhaps the most direct interpretation is that this value counts $(n-1,n)$-parking functions (parking functions with $n-1$ cars and $n$ spots) (this follows from the metered parking function interpretation by applying \cite[Proposition 2.4]{Daugherty}). By adding a $n$th car with preference $1$, $(n-1,n)$-parking functions can be viewed as classical $(n,n)$-parking functions with the sole restriction on the $n$th car. Utilizing the bijection of \cite{LRT-Weary}, we obtain the following result. 

\begin{cor} \label{cor:no-doublycovering}
    The number of rooted trees on $[n]_0$ with no doubly record-covering pairs is equal to the number of rooted trees on $[n]_0$ where $n$ is a child of $\circ$.
\end{cor}


\section{Record generating functions and the Cayley tree function}
\label{se:record_gf}

\bigskip
\subsection{Tree record generating functions}

 Fix a nonnegative integer $k$. 
Let $\mathcal{R}_k(z)$ be the generating function for rooted trees with precisely $k$ records:
\[
\mathcal{R}_k(z) = \sum_{n\ge 0} R(n,k) \ z^n
\]
and let $\mathcal{R}_{\ge k}(z)$ be the generating function for rooted trees at least $k$ records:
\[
\mathcal{R}_{\ge k}(z) = \sum_{l\ge k}\mathcal{R}_l(z) 
\]

\begin{cor} The generating function for trees with at least $k$ records is expressed in terms of the Cayley tree function as:
\[
\mathcal{R}_{\ge k}(z)  = \frac{\mathcal{T}^k(z)}{k}
\] 
\end{cor}

\begin{proof}
This follows immediately from Corollary \ref{pablos_lemma}, together with the cycle construction and its interpretation in terms of generating functions.
\end{proof}

In Proposition 3.8 of \cite{LRT-GF-records}, it was shown that the generating function for rooted trees with precisely $k$ records obeys the elegant equation
\begin{align}
\label{eq:k-dissymm}
\mathcal{R}_k = \frac{\mathcal{T}^k(z)}{k} - \frac{\mathcal{T}^{k+1}(z)}{k+1}
\end{align}
Indeed, for  $k=1$, Equation (\ref{eq:k-dissymm}) becomes the familiar expression for the generating function $\mathcal{U}(z)$ for unrooted trees:
\[
\mathcal{U}(z) =  \mathcal{T}(z) - \frac{\mathcal{T}^{2}(z)}{2}
\]
 (Observe that $\mathcal{U}(z)= \mathcal{R}_1(z)$ because
 trees rooted at their largest nodes are equinumerous with unrooted trees).  An equation given its “from-the-book” proof by Pierre Leroux's dissymmetry lemma \cite{Leroux_dissym, Leroux_dissym_2, gessel2023goodhuntingsproblemcounting}. Proposition 3.8 of \cite{LRT-GF-records} (Eq.~\ref{Eq:Leroux_generalization}) can be given a bijective proof that can be understood as a generalization of Leroux's dissymetry lemma

\begin{lem}[A generalization of Pierre Leroux's dissimmetry lemm]
\label{cor:gf_krecords}
\begin{align}
\label{Eq:Leroux_generalization}
 \frac{\mathcal{T}^k(z)}{k} =
 \mathcal{R}_k + \frac{\mathcal{T}^{k+1}(z)}{k+1} 
\end{align}
\end{lem}
\begin{proof}
There is a bijection between the set of endofunctions on $[n]$ of girth $k$ and the union of the set of trees of order $n$ with $k$ records and the set of endofunctions on $[n]$ of girth $k+1$.
Given an endofunction of girth $k$, the bijection of Corollary \ref{pablos_lemma} gives a tree $T$ with at least $k$ records. Either $T$ is in the set of trees of order $n$ with $k$ records, or it has at least $k+1$ records. In the latter case, the inverse of the same bijection gives an endofunction of girth $k+1$.

\end{proof}

\subsection{Differential equations and the record generating function}

In this subsection we provide a combinatorial interpretation for the derivative  $\mathcal{R}'_k(z)$ of the series 
$\mathcal{R}_k(z)$. This series satisfies
\[
\mathcal{R}'_k(z)=\frac{\mathcal{T}^k(z)}{z}.
\]
Indeed $\mathcal{R}'_1(z) = \frac{\mathcal{T}(z)}{z}$, as can be  checked directly by taking the derivative with respect to $z$ of $\mathcal{R}_1(z)$, the generating function for unrooted trees, and the expression of $\mathcal{R}_k(z)$ follows directly from Lemma \ref{cor:gf_krecords} by differentiating with respect to $z$:
\begin{align*}
\mathcal{R}'_k(z) &= \mathcal{T}^{k-1}(z) \mathcal{T}'(z) - \mathcal{T}^{k}(z) \mathcal{T}'(z) \\&= 
\mathcal{T}(z) \big( 
\mathcal{T}^{k-2}(z) \mathcal{T}'(z) - \mathcal{T}^{k-1}(z) \mathcal{T}'(z) \big)\\& = \mathcal{T}(z) \mathcal{R}_{k-1}'(z). 
\end{align*}

Our aim is thus to give a combinatorial proof of 
\[
\mathcal{R}'_k(z)=\mathcal{T}(z)\,\mathcal{R}'_{k-1}(z).
\]

 For this purpose, let us consider the class of rooted Cayley trees $\mathcal{R}^*_k$ with $k$ records, one of which is a marked virtual node, i.e. a node with no label and that does not contribute to the size, and such that the subtree of the virtual node does not contain any record or, in other terms, the path going from each record to the root does not cross the virtual node (except for the record corresponding to the virtual node), then we have the following lemma.
\begin{lem}
The series $\mathcal{R}'_k(z)$ is the generating function of the class $\mathcal{R}^*_k$
\label{lem:LemVirtNode}
\end{lem}
\begin{proof}There is a simple bijection between the trees of size $n-1$ in $\mathcal{R}^*_k$ and the trees of size $n$ in $\mathcal{R}_k$: this bijection  consists in giving the label $n$ to the virtual node of a tree of size $n-1$ from $\mathcal{R}^*_k$  and, conversely, in removing the label $n$ of a tree of size $n$ of $\mathcal{R}_k$ (see Figure~\ref{fig:BijCayley-virtual}). Since the virtual node is considered as a record and since it has no records in its sub-tree, the number of records does not change.
\begin{figure}[h!]
    \begin{center}
        \resizebox{0.75\textwidth}{!}{$%
        \vcenter{\hbox{\usetikzlibrary{positioning}
\tikzset{
  circ/.style={
    circle, draw, semithick, minimum size=1.5mm, inner sep=0pt, fill=black,
    }
}

\begin{tikzpicture}[
    arr/.style={semithick}
]
    \newcommand{\len}{1}

    \node[circ, label={right:$14$}] (s1) at (0,0) {};
    \node[circ, label={right:$2$}] (s11) at ($(s1) + (-135:\len)$) {};
    \node[circ, label={right:$5$}] (s12) at ($(s1) + (-90:\len)$) {};
    \node[circ, label={right:$6$}] (s13) at ($(s1) + (-45:\len)$) {};
    \node[circ, label={right:$8$}] (s111) at ($(s11) + (-135:\len)$) {};
    \node[circ, label={right:$11$}] (s112) at ($(s11) + (-90:\len)$) {};
    \node[circ, label={right:$13$}] (s113) at ($(s11) + (-45:\len)$) {};
    \node[circ, label={right:$7$}] (s1121) at ($(s112) + (-135:\len)$) {};
    \node[circ, label={right:$16$}] (s1122) at ($(s112) + (-45:\len)$) {};
    \node[circ, label={right:$9$}] (s1131) at ($(s113) + (-40:\len)$) {};
    \node[circ, label={right:$10$}] (s11311) at ($(s1131) + (-115:\len)$) {};
    \node[circ, label={right:$18$}] (s11312) at ($(s1131) + (-40:\len)$) {};         
     \node[circ, label={right:$15$}] (s113121) at ($(s11312) + (-135:\len)$) {};
     \node[circ, label={right:$17$}] (s113122) at ($(s11312) + (-90:\len)$) {};
    \node[circ, label={right:$1$}] (s1131221) at ($(s113122) + (-150:\len)$) {};
     \node[circ, label={right:$3$}] (s1131222) at ($(s113122) + (-120:\len)$) {};
      \node[circ, label={right:$4$}] (s1131223) at ($(s113122) + (-60:\len)$) {};
      \node[circ, label={right:$12$}] (s1131224) at ($(s113122) + (-30:\len)$) {};

    \draw[arr] (s1) -- (s11);
    \draw[arr] (s1) -- (s12);
    \draw[arr] (s1) -- (s13);
    \draw[arr] (s11) -- (s111);
    \draw[arr] (s11) -- (s112);
    \draw[arr] (s11) -- (s113);
    \draw[arr] (s112) -- (s1121);
    \draw[arr] (s112) -- (s1122);
    \draw[arr] (s113) -- (s1131);
    \draw[arr] (s1131) -- (s11311);
    \draw[arr] (s1131) -- (s11312);
    \draw[arr] (s11312) -- (s113121);
    \draw[arr] (s11312) -- (s113122);
    \draw[arr] (s113122) -- (s1131221);
    \draw[arr] (s113122) -- (s1131222);
    \draw[arr] (s113122) -- (s1131223);
    \draw[arr] (s113122) -- (s1131224);

\end{tikzpicture}}}%
        \hspace{1cm}%
        \longleftrightarrow%
        \hspace{1cm}%
        \vcenter{\hbox{\usetikzlibrary{positioning}
\tikzset{
  circ/.style={
    circle, draw, semithick, minimum size=1.5mm, inner sep=0pt, fill=black,
    }
}

\begin{tikzpicture}[
    arr/.style={semithick}
]
    \newcommand{\len}{1}

    \node[circ, label={right:$14$}] (s1) at (0,0) {};
    \node[circ, label={right:$2$}] (s11) at ($(s1) + (-135:\len)$) {};
    \node[circ, label={right:$5$}] (s12) at ($(s1) + (-90:\len)$) {};
    \node[circ, label={right:$6$}] (s13) at ($(s1) + (-45:\len)$) {};
    \node[circ, label={right:$8$}] (s111) at ($(s11) + (-135:\len)$) {};
    \node[circ, label={right:$11$}] (s112) at ($(s11) + (-90:\len)$) {};
    \node[circ, label={right:$13$}] (s113) at ($(s11) + (-40:\len)$) {};
    \node[circ, label={right:$7$}] (s1121) at ($(s112) + (-135:\len)$) {};
    \node[circ, label={right:$18$}] (s1122) at ($(s112) + (-45:\len)$) {};
    \node[circ, label={right:$16$}] (s11221) at ($(s1122) + (-135:\len)$) {};
    \node[circ, label={right:$9$}] (s1131) at ($(s113) + (-45:\len)$) {};
    \node[circ, label={right:$10$}] (s11311) at ($(s1131) + (-115:\len)$) {};
    \node[circ, fill=white, draw=black] (s11312) at ($(s1131) + (-45:\len)$) {};
     \node[circ, label={right:$15$}] (s113121) at ($(s11312) + (-135:\len)$) {};
     \node[circ, label={right:$17$}] (s113122) at ($(s11312) + (-90:\len)$) {};
     \node[circ, label={right:$1$}] (s1131221) at ($(s113122) + (-150:\len)$) {};
     \node[circ, label={right:$3$}] (s1131222) at ($(s113122) + (-120:\len)$) {};
      \node[circ, label={right:$4$}] (s1131223) at ($(s113122) + (-60:\len)$) {};
      \node[circ, label={right:$12$}] (s1131224) at ($(s113122) + (-30:\len)$) {};

    \draw[arr] (s1) -- (s11);
    \draw[arr] (s1) -- (s12);
    \draw[arr] (s1) -- (s13);
    \draw[arr] (s11) -- (s111);
    \draw[arr] (s11) -- (s112);
    \draw[arr] (s11) -- (s113);
    \draw[arr] (s112) -- (s1121);
    \draw[arr] (s112) -- (s1122);
    \draw[arr] (s1122) -- (s11221);
    \draw[arr] (s113) -- (s1131);
    \draw[arr] (s1131) -- (s11311);
    \draw[arr] (s1131) -- (s11312);
    \draw[arr] (s11312) -- (s113121);
    \draw[arr] (s11312) -- (s113122);
    \draw[arr] (s113122) -- (s1131221);
    \draw[arr] (s113122) -- (s1131222);
    \draw[arr] (s113122) -- (s1131223);
    \draw[arr] (s113122) -- (s1131224);

\end{tikzpicture}}}
        $} 
    \end{center}
  \caption{An example of the bijection between Cayley trees with $n$ labels and Cayley trees with $n-1$ labels and a virtual node. The virtual node is in white.}
    \label{fig:BijCayley-virtual}
\end{figure}
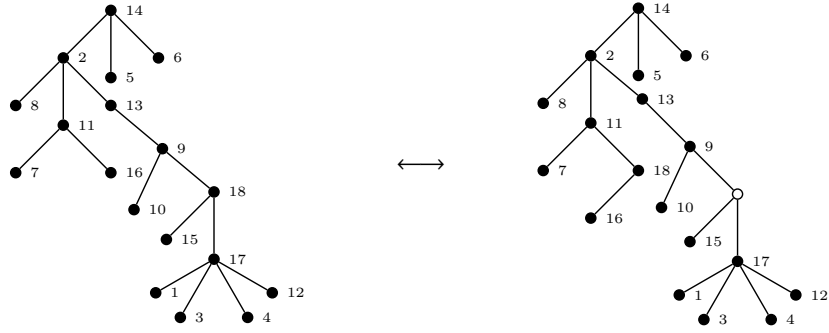

Then, by using the bijection we have the following:
\[
\sum_{T \in \mathcal{R}^*_k} \frac{z^{|T|}}{|T|!}=\sum_{T' \in \mathcal{R}_k} \frac{z^{|T'|-1}}{(|T'|-1)!}=\mathcal{R}'_k(z).
\]
\end{proof}

\begin{figure}[h!]
    \centering
\resizebox{1\textwidth}{!}{$%
  \vcenter{\hbox{\usetikzlibrary{positioning}
\tikzset{
  circ/.style={
    circle, draw, semithick, minimum size=1.5mm, inner sep=0pt, fill=black,
    }
}

\begin{tikzpicture}[
    arr/.style={semithick}
]
    \newcommand{\len}{1}

    \node[circ, label={right:$14$}] (s1) at (0,0) {};
    \node[circ, label={right:$2$}] (s11) at ($(s1) + (-135:\len)$) {};
    \node[circ, label={right:$5$}] (s12) at ($(s1) + (-90:\len)$) {};
    \node[circ, label={right:$6$}] (s13) at ($(s1) + (-45:\len)$) {};
    \node[circ, label={right:$8$}] (s111) at ($(s11) + (-135:\len)$) {};
    \node[circ, label={right:$11$}] (s112) at ($(s11) + (-90:\len)$) {};
    \node[circ, label={right:$13$}] (s113) at ($(s11) + (-40:\len)$) {};
    \node[circ, label={right:$7$}] (s1121) at ($(s112) + (-135:\len)$) {};
    \node[circ, label={right:$18$}] (s1122) at ($(s112) + (-45:\len)$) {};
    \node[circ, label={right:$16$}] (s11221) at ($(s1122) + (-135:\len)$) {};
    \node[circ, label={right:$9$}] (s1131) at ($(s113) + (-45:\len)$) {};
    \node[circ, label={right:$10$}] (s11311) at ($(s1131) + (-115:\len)$) {};
    \node[circ, fill=white, draw=black] (s11312) at ($(s1131) + (-45:\len)$) {};
     \node[circ, label={right:$15$}] (s113121) at ($(s11312) + (-135:\len)$) {};
     \node[circ, label={right:$17$}] (s113122) at ($(s11312) + (-90:\len)$) {};
     \node[circ, label={right:$1$}] (s1131221) at ($(s113122) + (-150:\len)$) {};
     \node[circ, label={right:$3$}] (s1131222) at ($(s113122) + (-120:\len)$) {};
      \node[circ, label={right:$4$}] (s1131223) at ($(s113122) + (-60:\len)$) {};
      \node[circ, label={right:$12$}] (s1131224) at ($(s113122) + (-30:\len)$) {};

    \draw[arr] (s1) -- (s11);
    \draw[arr] (s1) -- (s12);
    \draw[arr] (s1) -- (s13);
    \draw[arr] (s11) -- (s111);
    \draw[arr] (s11) -- (s112);
    \draw[arr] (s11) -- (s113);
    \draw[arr] (s112) -- (s1121);
    \draw[arr] (s112) -- (s1122);
    \draw[arr] (s1122) -- (s11221);
    \draw[arr] (s113) -- (s1131);
    \draw[arr] (s1131) -- (s11311);
    \draw[arr] (s1131) -- (s11312);
    \draw[arr] (s11312) -- (s113121);
    \draw[arr] (s11312) -- (s113122);
    \draw[arr] (s113122) -- (s1131221);
    \draw[arr] (s113122) -- (s1131222);
    \draw[arr] (s113122) -- (s1131223);
    \draw[arr] (s113122) -- (s1131224);

\end{tikzpicture}}}%
  \hspace{1cm}%
  \longleftrightarrow%
  \hspace{1cm}%
  \vcenter{\hbox{\usetikzlibrary{positioning}
\tikzset{
  circ/.style={
    circle, draw, semithick, minimum size=1.5mm, inner sep=0pt, fill=black,
    }
}

\begin{tikzpicture}[
    arr/.style={semithick}
]
    \newcommand{\len}{1}

    \node[circ, label={right:$14$}] (s1) at (0,0) {};
    \node[circ, label={right:$2$}] (s11) at ($(s1) + (-135:\len)$) {};
    \node[circ, label={right:$5$}] (s12) at ($(s1) + (-90:\len)$) {};
    \node[circ, label={right:$6$}] (s13) at ($(s1) + (-45:\len)$) {};
    \node[circ, label={right:$8$}] (s111) at ($(s11) + (-135:\len)$) {};
    \node[circ, label={right:$11$}] (s112) at ($(s11) + (-90:\len)$) {};
    \node[circ, label={right:$13$}] (s113) at ($(s11) + (-40:\len)$) {};
    \node[circ, label={right:$7$}] (s1121) at ($(s112) + (-135:\len)$) {};
    \node[circ, label={right:$18$}] (s1122) at ($(s112) + (-45:\len)$) {};
    \node[circ, label={right:$16$}] (s11221) at ($(s1122) + (-135:\len)$) {};
    \node[circ, fill=white, label={right:$9$}] (s1131) at ($(s113) + (-45:\len)$) {};
    \node[circ, label={right:$10$}] (s11311) at ($(s1131) + (-115:\len)$) {};
    \node[circ, fill=white, draw=blue] (s11312) at ($(s1131) + (-45:\len)$) {};
     \node[circ, fill=blue, draw=blue, label={right:$\textcolor{blue}{15}$}] (s113121) at ($(s11312) + (-135:\len)$) {};
     \node[circ, fill=blue, draw=blue,label={right:$\textcolor{blue}{17}$}] (s113122) at ($(s11312) + (-90:\len)$) {};
     \node[circ, fill=blue, draw=blue,label={right:$\textcolor{blue}{1}$}] (s1131221) at ($(s113122) + (-150:\len)$) {};
     \node[circ, fill=blue, draw=blue,label={right:$\textcolor{blue}{3}$}] (s1131222) at ($(s113122) + (-120:\len)$) {};
      \node[circ, fill=blue, draw=blue,label={right:$\textcolor{blue}{4}$}] (s1131223) at ($(s113122) + (-60:\len)$) {};
      \node[circ, fill=blue, draw=blue,label={right:$\textcolor{blue}{12}$}] (s1131224) at ($(s113122) + (-30:\len)$) {};

    \draw[arr] (s1) -- (s11);
    \draw[arr] (s1) -- (s12);
    \draw[arr] (s1) -- (s13);
    \draw[arr] (s11) -- (s111);
    \draw[arr] (s11) -- (s112);
    \draw[arr] (s11) -- (s113);
    \draw[arr] (s112) -- (s1121);
    \draw[arr] (s112) -- (s1122);
    \draw[arr] (s1122) -- (s11221);
    \draw[arr] (s113) -- (s1131);
    \draw[arr] (s1131) -- (s11311);
    
    \draw[arr, draw=blue] (s11312) -- (s113121);
    \draw[arr, draw=blue] (s11312) -- (s113122);
    \draw[arr, draw=blue] (s113122) -- (s1131221);
    \draw[arr, draw=blue] (s113122) -- (s1131222);
    \draw[arr, draw=blue] (s113122) -- (s1131223);
    \draw[arr, draw=blue] (s113122) -- (s1131224);

\end{tikzpicture}}}%
   \hspace{1cm}%
  \longleftrightarrow%
  \hspace{1cm}%
  \vcenter{\hbox{\usetikzlibrary{positioning}
\tikzset{
  circ/.style={
    circle, draw, semithick, minimum size=1.5mm, inner sep=0pt, fill=black,
    }
}

\begin{tikzpicture}[
    arr/.style={semithick}
]
    \newcommand{\len}{1}

    \node[circ, label={right:$16$}] (s1) at (0,0) {};
    \node[circ, label={right:$2$}] (s11) at ($(s1) + (-135:\len)$) {};
    \node[circ, label={right:$5$}] (s12) at ($(s1) + (-90:\len)$) {};
    \node[circ, label={right:$6$}] (s13) at ($(s1) + (-45:\len)$) {};
    \node[circ, label={right:$8$}] (s111) at ($(s11) + (-135:\len)$) {};
    \node[circ, label={right:$13$}] (s112) at ($(s11) + (-90:\len)$) {};
    \node[circ, label={right:$14$}] (s113) at ($(s11) + (-40:\len)$) {};
    \node[circ, label={right:$7$}] (s1121) at ($(s112) + (-135:\len)$) {};
    \node[circ, fill=white, draw=black] (s1122) at ($(s112) + (-45:\len)$) {};
    \node[circ, label={right:$18$}] (s11221) at ($(s1122) + (-135:\len)$) {};
    \node[circ, label={right:$10$}] (s1131) at ($(s113) + (-45:\len)$) {};
    \node[circ, label={right:$11$}] (s11311) at ($(s1131) + (-115:\len)$) {};
    \node[circ, fill=blue, draw=blue, label={right:\textcolor{blue}{$9$}}] (s11312) at ($(s1131) + (-45:\len)$) {};
     \node[circ, fill=blue, draw=blue, label={right:\textcolor{blue}{$15$}}] (s113121) at ($(s11312) + (-135:\len)$) {};
     \node[circ, fill=blue, draw=blue,, label={right:\textcolor{blue}{$17$}}] (s113122) at ($(s11312) + (-90:\len)$) {};
     \node[circ, fill=blue, draw=blue,, label={right:\textcolor{blue}{$1$}}] (s1131221) at ($(s113122) + (-150:\len)$) {};
     \node[circ, fill=blue, draw=blue, label={right:\textcolor{blue}{$3$}}] (s1131222) at ($(s113122) + (-120:\len)$) {};
      \node[circ, fill=blue, draw=blue, label={right:\textcolor{blue}{$4$}}] (s1131223) at ($(s113122) + (-60:\len)$) {};
      \node[circ, fill=blue, draw=blue, label={right:\textcolor{blue}{$12$}}] (s1131224) at ($(s113122) + (-30:\len)$) {};

    \draw[arr] (s1) -- (s11);
    \draw[arr] (s1) -- (s12);
    \draw[arr] (s1) -- (s13);
    \draw[arr] (s11) -- (s111);
    \draw[arr] (s11) -- (s112);
    \draw[arr] (s11) -- (s113);
    \draw[arr] (s112) -- (s1121);
    \draw[arr] (s112) -- (s1122);
    \draw[arr] (s1122) -- (s11221);
    \draw[arr] (s113) -- (s1131);
    \draw[arr] (s1131) -- (s11311);
    \draw[arr, draw=blue] (s11312) -- (s113121);
    \draw[arr, draw=blue] (s11312) -- (s113122);
    \draw[arr, draw=blue] (s113122) -- (s1131221);
    \draw[arr, draw=blue] (s113122) -- (s1131222);
    \draw[arr, draw=blue] (s113122) -- (s1131223);
    \draw[arr, draw=blue] (s113122) -- (s1131224);

\end{tikzpicture}}}%
  $}
  \caption{The decomposition $\phi$ of a Cayley tree with a virtual node. The trees obtained are $\bar S$, in black, and $\bar S'$, in blue. In the center there are the intermediary trees $S$ and $S'$.}
  \label{fig:DecCayleyVirtual}
\end{figure}
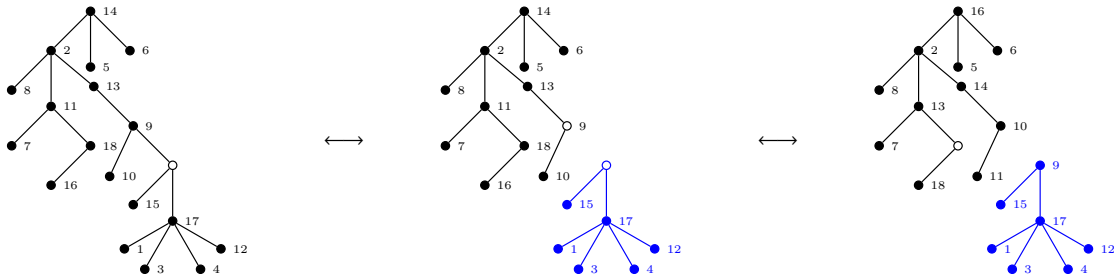

Let us now define an application $\phi$ mapping each element of $\mathcal{R}^*_k$ to a pair consisting after standardization in a Cayley tree and an element of $\mathcal{R}^*_{k-1}$:
Given a tree $T$ belonging to $\mathcal{R}^*_k$,  we mark the parent of the virtual node and cut the tree into two sub-trees at the level of the virtual node. 
We denote by  $S'$  the sub-tree of $T$ whose root is the virtual node, and  by  $S$ the sub-tree obtained from $T$ by stripping the sub-tree $S'$ and marking the node where it was attached (see Figure~\ref{fig:DecCayleyVirtual}, left-hand side). Then we transfer the label $\ell$ of the marked node of $S$ to the root of $S'$ to form a labeled tree $\bar S'$. Finally we shift all the labels $\ell_1<\ldots<\ell_p$ of $S$ that are greater than $\ell$ and remove the mark to get a tree $\bar S$: the marked node of $S$ gets label $\ell_i$ in $\bar S$, the node with label $\ell_i$  gets label $\ell_{i+1}$, and the node with maximal label $\ell_p$ in $S$ becomes a virtual (unlabeled) node in $\bar S$. 
\begin{p}
       The mapping $\phi$ that sends $T$ to $(\bar S,\bar S')$ is a bijection between
    \begin{itemize}
        \item trees of $\mathcal{R}^*_k$ of size $n$,
        \item and pairs of trees $(\bar S,\bar S')$ with total label set $\{1,\ldots,n\}$ such that the standardization of $\bar S$ belongs to $\mathcal{R}^*_k$ and the standardization of $\bar S'$ is a Cayley tree.
    \end{itemize}
\end{p}
\begin{proof}
Observe that the tree $S$ has $k-1$ records since by definition of the class  $\mathcal{R}^*_k$ there are no records in the sub-tree starting at the virtual node, except for the virtual node itself, and the shift of labels preserves the records: in particular the virtual node of $\bar S$ remains a record since it takes the place of the maximum of $S$.

Clearly $S$ and $S'$ can be recovered from $\bar S$ and $\bar S'$ upon removing label $\ell$ at the root of $\bar S'$ and inserting it at the first larger label $\ell_1>\ell$ in $\bar S$ and shifting back the larger labels $\ell_1<\ldots<\ell_p$ until $\ell_p$ is put back on the virtual node of $\bar S$. Then $S$ and $S'$ are glued together to recover $T$. This proves that the mapping $\phi$ is injective. 

Conversely given an pair of trees $(\bar S,\bar S')$ as in the proposition the inverse construction yields a tree $T$ with a virtual node and $k$ records that are not in the sub-tree of the virtual node, that is a tree of $\mathcal{R}^*_k$, and $\phi(T)=(\bar S,\bar S')$.
\end{proof}

At this point, the fundamental theorem of Calculus, implies that an integral form for $\mathcal{R}_k$ in terms of Cayley's tree function.

\begin{cor}[Proposition 3.4 of \cite{LRT-GF-records}]
\label{le:trees_k_records}
The generating function for   trees with precisely $k$ records can be expressed in terms of the Cayley tree function as
\begin{align*}
\mathcal{R}_k = 
\int_0^z  \frac{1}{s}
  \mathcal{T}^k(s)  \ ds
  \end{align*}
\end{cor}

\subsection{Forest record generating function}
\label{subsec:forest record generating function}

We say that a rooted forest has been \emph{marked} when a record has been distinguished at each  of its constituent trees. 
\begin{cor}[Theorem \ref{thm:connected endofunctions}]
There is a bijection between the number of marked forests of order $n$ and the set of endofunctions of $[n]$, where the height of the distinguished records determines the girth of the connected components of the endofunction. 
\end{cor}

\begin{proof} Let $F$ be a marked forest.
If we apply the bijection of Theorem \ref{thm:connected endofunctions} to each constituent tree of $F$, the result is an endofunction $f$ with the same number of connected components as $F$. Finally,  Theorem \ref{thm:connected endofunctions}  also tells us how to recover the height of the selected record nodes of $F$ from the girth of the connected components of $f$.

 \end{proof}

Define the \emph{multiplicative records} of a rooted forest $F = \{T_1, T_2, \ldots, T_k\}$ as the product of the number of records of its connected components:
$
\text{mrec}(F) =\rec(T_1) \ldots \rec(T_k)$.

 \begin{cor}
The total number of marked forests of order $n$ is $n^n.$  Moreover,
\begin{align}
\label{Eq:record_forests}
\sum_F \mrec(F) = n^n
\end{align}
where the sum is taken over all rooted forests $F$ with labels $[n]$.
\end{cor}

We can do better and keep track of the girth of the connected components. Recall that the type on an endofunction  $f$ is the partition whose parts are the girths of its connected components.
Let
\begin{align*}
C(z,t) = 
\sum_{(T,r) }  \frac{z^{|T|} t^{ht(r)-1}}{|T|! },
& \qquad \qquad 
F(z,t) = 
\sum_{(F, \vec{r}) }  \frac{z^{|F|} t^{ht(\vec{r})-\vec{1}}}{|F|!}
\end{align*}
where the first sum is taken over the set of rooted trees with a selected node, and the second sum over all  rooted forests with a selected record in each constituent tree (that we order according to the labels of their roots). Theorem \ref{thm:connected endofunctions} implies that
\[C(z, t) = \log\left(\frac{1}{1 - t  \mathcal{T}(z)}\right)
\] is the well-studied exponential generating function for connected endofunctions. For more information, see \url{https://oeis.org/A201685}.
On the other hand, since the forest record numbers are multiplicative (Eq.~(\ref{Eq:record_forests})) the exponential formula implies that:
\[
    F(z, t)=  \frac{1}{1 - t\mathcal{T}(z)}
\]
 The reader can find more information on these generating functions, including some asymptotic results in  page 130 of \cite{Flajolet}.

On the other hand, recall that the \emph{distance} between any pair of nodes $v_1,v_2$ is the length of the unique path between $v_1$ and $v_2$. Therefore, $F(z, 1)^2/2$ is the exponential generating function for the total sum of the distances between all nodes of all unrooted trees of order $n$:
\[
\frac{F(z,1)^2}{2}
=
\sum_{n\ge 1} 
\sum_{T}  
d(T)
\frac{z^n}{n!} 
\]
where the second  sum is taken over all unrooted labeled trees $T$ on $n$ nodes, and $d(T)$ is defined as the sum of the distances between all pairs of nodes of $T$. 
That is, 
\[
d(T) = \sum_{1\le i<j \le n} d(i, j).
\]
This is sequence \url{https://oeis.org/A036276}.

There is a simple bijective proof of this result. Let $T$ be a  tree of order $n.$ 
 Choose a two vertices $i$ and $j$ in $T$, and an edge between them. This can be done in $\sum_{1\le i<j \le n} d(i, j)$ different ways (if $i=j$ this can be done in zero ways, hence, it is pointless to ask that $i\ne j$). 
 On the other hand, cutting at  the selected edge $e$ splits $T$ into two trees $T_i$ and $T_j$, with $T_i$  rooted at $i$ and with the endpoint of $e$ in this component selected, and $T_j$ rooted at $j$, and with the node at the other endpoint of $e$ selected. 
The number of trees with two selected nodes (that may coincide) is $n^n$.

 Let $E(z)$ denote the generating function for endofunctions from an $n$-set to itself. Then,
\[
E(z)= \sum_{n\ge 1} n^n \frac{z^n}{n!}
= \sum_{n\ge 1} n^{n-1} \frac{z^n}{(n-1)! }
= z \sum_{n\ge 1} n^{n-1} \frac{z^{n-1}}{(n-1)!} 
= z \ \partial_z \mathcal{T}(z)
\]
We claim that $E(z)^2/2$ is the gf for
\url{https://oeis.org/A036276}

\[
\frac{E(z)^2}{2}
=
\sum_{n\ge 1} 
\sum_{T}  
d(T)
\frac{z^n}{n!} 
\]
where the second  sum is taken over all  $n^{n-2}$ unrooted labeled trees $T$ on $n$ nodes. 
 Note that $E(z)^2$ coincides with $W(z)$, the generating function for rooted trees with a catalyst (Lemma \ref{lem:Riordan_Sloane}). The equality above is therefore equivalent to the following statement: the sum of the total heights over all rooted trees is exactly twice the sum of the total distance of all unrooted trees.
There is a simple combinatorial proof of this result.  Let $T$ be an unrooted tree. The total distance $d(T)$ equals the sum of the lengths of all paths in $T$. 
On the other hand, when summing the total height over all rooted versions of $T$, the length of each path in $T$ contributes exactly twice, as either of its endpoints may serve as the root.


\section{Coda. Cayley's  forest  formula.}
\label{App:Cayley_forest}
To close, we propose a new proof of Cayley’s forest formula using the record decomposition.
Given any rooted forest \(F\) with $n$ nodes, we construct a rooted tree $T$ by attaching a new vertex \(\circ\) to all roots of \(F\), and labelling $\circ$ with 0. We stress that  $0$ can only appear as the label of the root  and that $\circ$ is never considered to be a record. We then  say that $T$ is labelled with $[n]_0$.

The roots of $F$ are referred to as \emph{children of the root} of $T$.  This well-known construction  establishes a bijection between forests (with $n$ nodes) rooted at $[k]$ and rooted trees (with $n$ plus one node) having $[k]$ as children of the root

The following lemma provides an elegant characterisation for the set of rooted trees with  nodes $[n]_0$ and  \([k]\) as set of children of the root in terms of the lexicographic order of its record codes.

\begin{lem}
    \label{lem:blob-iff}
Let $T$ be a rooted tree labelled with nodes $[n]_0$. Then, the set of children of the root of  $T$ is $[k]$ if and only if 
  \[\circ^{k-1}1^{n - k} \leq \blob(T) \leq  \circ^{k-1}kn^{n-k-1}\]
   where $\leq$ denotes the lexicographic ordering.
\end{lem}

\begin{proof} Let $T$ be a   tree labelled with $[n]_0$, and let 
$[k]$ be the set of children of the root.
    Write $\blob(T) = x_1 x_2\cdots x_{n-1}$  for  the second line of the record code of $T$.
    By hypothesis, the elements of $[k]$ are children of the root of $T$, thus $x_i = f_T(i + 1) = \circ$ for $i \in [k - 1]$. 
    On the other hand, since $k$ is a record, it must happen that $x_k \leq k$. 
     Lastly, since there are $k$ children of the root, the multiplicity of $\circ$ in $C$ must be $k - 1$, the degree of $\circ$ in $T$ minus one. Hence $x_i \geq 1$ for all $i \in [k + 1, n - 1]$.   
    
    Conversely, assume that $C = x_1x_2\cdots x_{n-1}$ is a code satisfying that $ \circ^{k-1} 1^{n-k}\leq C \leq \circ^{k-1}kn^{n-k-1}$, and let $T$ be the tree encoded by $C$, that is, such that $\blob(T) = C$. Then, all vertices in $[k]$ are records of $T$. Thus $f_T(1) = \circ$ (since $1$ is the smallest record of $T$) and  $f_T(i) = x_{i - 1} = \circ$ for $i \in [2, k]$. Since the degree of $\circ$ in $T$ is $k$ and we have already determined its $k$ children, the parent of any other vertex cannot be $\circ$. Thus, $f_T(i) \neq \circ$ for $i \in [k + 1, n]$.
\end{proof}

 Cayley's  forest  formula immediately follows after observing that the number of record codes that satisfy the conditions stated in Lemma \ref{lem:blob-iff} is $k n^{n-k-1}$.

\begin{cor}[Cayley's  forest  formula, \cite{Cayley}]
    \label{prop: bijection [k]-rooted}
    Let  $k \leq n$ be natural numbers. 
        The number of forests on $n$ vertices rooted at $[k]$ is
    \begin{equation}
    \label{eqn:|B(n,k)|}
     k n^{n-k-1}.
    \end{equation}
\end{cor}

We remark that having a closed formula for these numbers, it is routine to verify that for each $n$ the sequence $k n^{n-k-1}$ is log-concave. \\

Since the number of forest on $n$ nodes and rooted at $S\subseteq[n]$ only depends on the cardinality of $S$, we easily recover this form of Cayley's  forest  formula.
\begin{cor}[Cayley, \cite{Cayley}] The number of spanning forests on $n$ vertices and $k$ trees is:
\[
{n-1 \choose k-1} n^{n-k} 
\]
\end{cor}

\section{Final Comments}

\begin{enumerate}[wide, labelindent=0pt]

\item
We pose as a question of interest the construction of a direct bijective proof of Corollary~\ref{cor:no-doublycovering}.

\item In Section \ref{subsec:forest record generating function} we study the sum of pairwise distances across all unrooted labelled trees on $n$ vertices. Remarkably, this quantity appears to coincide with the total number of defects in all defective parking functions of length $n$, see \url{https://oeis.org/A036276}.
On the other hand, the number of trees with no doubly record-covering pairs seems to coincide with the number of metered parking functions, see Section \ref{sec:doublycovering-records} and \url{https://oeis.org/A007334}. 
A direct combinatorial proof of these facts would be highly desirable.

\end{enumerate}


\section*{Acknowledgements}
We gratefully acknowledge the support and engaging discussions with Cyril Bandelier, and the insightful remarks of Andrea Sportiello.

This work has been partially supported by Grants PID2020-117843GB-
I00 and PID2024-157173NB-I00 funded by MCIN/AEI/10.13039/501100011033 and by FEDER, UE. Author ST is supported by \href{10.3030/101070558}{FoQaCiA} which is funded by the European Union and NSERC.

\bibliographystyle{alpha}
\bibliography{references.bib}

\end{document}